\documentclass[]{interact}

\usepackage{epstopdf}% To incorporate .eps illustrations using PDFLaTeX, etc.
\usepackage{subfigure}% Support for small, `sub' figures and tables

\usepackage{natbib}% Citation support using natbib.sty
\bibpunct[, ]{(}{)}{;}{a}{}{,}% Citation support using natbib.sty
\renewcommand\bibfont{\fontsize{10}{12}\selectfont}% Bibliography support using natbib.sty

\theoremstyle{plain}% Theorem-like structures provided by amsthm.sty

\theoremstyle{definition}

\theoremstyle{remark}

\usepackage[T1]{fontenc} 
\usepackage[utf8]{inputenc}
\usepackage[english]{babel}
\usepackage{amsmath,amssymb,amsfonts}
\usepackage{algorithm}
\usepackage{algorithmicx}
\usepackage{graphicx}
\usepackage{textcomp}
\usepackage{xcolor}
\usepackage{physics}
\usepackage{booktabs}
\usepackage{threeparttable}
\usepackage{multirow}
\usepackage{longtable}
\usepackage{setspace}
\usepackage{placeins}
\usepackage{lipsum}
\usepackage{pdfpages}
\usepackage{hyperref}
\usepackage{color,soul}
\usepackage[normalem]{ulem}
\usepackage[version=4]{mhchem}
\usepackage[euler]{textgreek}
\usepackage{titlesec}
\renewcommand{\vec}{\boldsymbol}

\begin{document}

\articletype{ARTICLE TEMPLATE}

\title{Nonlinear Programming Solvers for Unconstrained and Constrained Optimization Problems: a Benchmark Analysis}

\author{
\name{
Giovanni Lavezzi\textsuperscript{a}
\thanks{CONTACT: Giovanni Lavezzi, Ph.D. student. Email: giovanni.lavezzi@sdstate.edu} 
, Kidus Guye\textsuperscript{b} 
\thanks{CONTACT: Kidus Guye, Ph.D. student. Email: g.kidus@wustl.edu} 
and Marco Ciarci\`a\textsuperscript{a}
\thanks{Corresponding author: Marco Ciarci\`a, Assistant professor. Email: marco.ciarcia@sdstate.edu} 
}
\affil{
\textsuperscript{a}Department of Mechanical Engineering, South Dakota State University, 1451 Stadium Rd, Brookings, South Dakota, 57006, USA; 
\textsuperscript{b}Department of Mechanical Engineering and Materials Science, Washington University in St. Louis, 1 Brookings Dr, St. Louis, 63130, Missouri, USA
}
}

\maketitle

% \title[NLP Solvers for Unconstrained and Constrained Problems]{Nonlinear Programming Solvers for Unconstrained and Constrained Optimization Problems: a Benchmark Analysis}

% \author*[1]{\fnm{Giovanni} \sur{Lavezzi}}\email{giovanni.lavezzi@sdstate.edu}

% \author[2]{\fnm{Kidus} \sur{Guye}}\email{g.kidus@wustl.edu}

% \author[3]{\fnm{Marco} \sur{Ciarci\`a}}\email{marco.ciarcia@sdstate.edu}

% \affil*[1]{\tanm{Ph.D. student}, \orgdiv{Department of Mechanical Engineering}, \orgname{South Dakota State University}, \orgaddress{\street{1451 Stadium Rd}, \city{Brookings}, \postcode{57006}, \state{South Dakota}, \country{USA}}}

% \affil[2]{\tanm{Ph.D. student}, \orgdiv{Department of Mechanical Engineering and Materials Science}, \orgname{Washington University in St. Louis}, \orgaddress{\street{1 Brookings Dr}, \city{St. Louis}, \postcode{63130}, \state{Missouri}, \country{USA}}}

% \affil[3]{\tanm{Assistant professor}, \orgdiv{Department of Mechanical Engineering}, \orgname{South Dakota State University}, \orgaddress{\street{1451 Stadium Rd}, \city{Brookings}, \postcode{57006}, \state{South Dakota}, \country{USA}}}


\begin{abstract}
In this paper we propose a set of guidelines to select a solver for the solution of nonlinear programming problems. With this in mind, we present a comparison of the convergence performances of commonly used solvers for both unconstrained and constrained nonlinear programming problems. The comparison involves accuracy, convergence rate, and convergence speed. Because of its popularity among research teams in academia and industry, MATLAB is used as common implementation platform for the solvers. Our study includes solvers which are either freely available, or require a license, or are fully described in literature. In addition, we differentiate solvers if they allow the selection of different optimal search methods. As result, we examine the performances of 23 algorithms to solve 60 benchmark problems. To enrich our analysis, we will describe how, and to what extent, convergence speed and accuracy can be improved by changing the inner settings of each solver.
\end{abstract}


\section{Introduction}
The current technological era prioritizes, more than ever, high performance and efficiency of complex processes controlled by a set of variables. Examples of these processes are \citep{Lasdon1980,Ignacio1996,Charalambous1979,Grossmann1997,Wu1992,Wansuo2010,Rustagi1994,Ziemba1975}: engineering designs, chemical plant reactions, manufacturing processes, grid power management, power generation/conversion process, path planning for autonomous vehicles, climate simulations, etc. Quite often, the search for the best performance, or the highest efficiency, can be transcribed into the form of a Nonlinear Programming (NLP) problem. Namely, the need to minimize (or maximize) a scalar cost function subjected to a set of constraints. In some instances these functions are linear but, in general, one or both of them are characterized by nonlinearities. For simple, one-time use problems, one might successfully use any of the solver available, like \textit{fmincon} in MATLAB \citep{Fmincon,matlab}. Nevertheless, if the NLP derives from some specific applications, like real-time process optimization, then the solver choice begs a more accurate selection. 

The first research efforts toward the characterization of optimization solvers started in the 60s'. In \citep{Box1966} the authors compare eight solvers on twenty benchmark unconstrained NLP problems containing up to 20 variables. Notably, they illustrate techniques to transform particular constrained NLP into equivalent unconstrained problems. The authors of \citep{Levy1976} analyze the convergence properties of two gradient-based solvers applied to 16 test problems. In the last few decades, with the development of new methodologies and optimization applications, more studies aimed to illustrate difference in performance among NLP solvers. Schittkowski et.al. \citep{Schit1994} performed the comparison of eleven different mathematical programming codes applied to structural optimization through finite element analysis. George et. al. summarize a qualitative comparison of few optimization methodologies reported by several other sources \citep{George1969}.  In the research document prepared by Sandia National Laboratory \citep{Gearhart2013}, a study was conducted on four open source Linear Programming (LP) solvers applied to 201 benchmark problems. In \citep{Kronqvist2018} Kronqvist et.al. carried out a performance comparison of mixed integer NLP solvers limited to convex benchmark problems.  Pratiksha Saxena presents a comparison between linear and nonlinear programming techniques on the diet formulation for animals \citep{Saxena2011}. Another work by Hannes Pucher uses the programming language R to analyze multiple nonlinear optimization methods applied to real life problems \citep{Pucher}. State-of-the-art optimization methods were used to compare on their application on L1-regularized classifiers \citep{Yuan2010}. On a similar note, multiple global optimization solvers were compared on a work done by Arnold Neumaier \citep{Neumaier2005}. Authors from \citep{Oba1997,Ilhagga1996,Haupt1995,Hamdy2016} conducted a performance comparison of optimization techniques for specific applications which includes aerodynamic shape design, integrated manufacturing planning and scheduling, solving electromagnetic problems, and building energy design problems, respectively. Similarly, Frank et. al. conducted a comparison between three optimization methods for solving aerodynamic design problems \citep{Frank1992}. 
In \citep{Karaboga2007}, Karaboga et. al. compares the performance of artificial bee colony algorithm with the differential evolution, evolutionary and particle swarm optimization algorithms using multi-dimensional numerical problems.

In this paper we want to provide an explicit comparison of a set of NLP solvers. We include in our comparison popular solvers readily available in MATLAB, a few gradient descent methods that have been extensively used in literature, and a particle swarm optimization. Because of its widespread use among research groups, both in academia and private sector, we have decided to use MATLAB as common implementation platform. For this reason we will focus on all the solvers that are either written on or can be implemented in MATLAB. The NLP problems used in this comparison have been selected amongst the standard benchmark problems \citep{testcollect,benchmark,handbook} with up to thirty variables and a up to nine scalar constraints. The paper is organized as follows. Section \ref{section:s2} describes the statement of unconstrained and contrained NLP problems. In Section \ref{section:s3}, we enumerate the NLP solvers included in our analysis and their main features. Subsequently, an overview of the different convergence metrics, and the solvers implementations is carried out in Section \ref{section:s4}. The results of the comparison with the benchmark equations is discussed in Section \ref{section:s5}. Finally, the main contributions of the paper are outlined in Section \ref{section:s6}.

\section{Nonlinear programming problem statements}
\label{section:s2}
In general, a constrained NLP problem aims to minimize a nonlinear real scalar objective function, with respect a set of variables, while satisfying a set of nonlinear constraints. If the problem entails the minimization of a function without the presence of constraints the problem is defined as unconstrained \citep{Wright2006}.
In the following section, the general form of a nonlinear unconstrained and constrained optimization problems in the minimization format form are thoroughly stated.

\subsection{Unconstrained optimization problem}
\subsubsection{Statement}
Let $\vec{x}\in\mathbb{R}^{n}$ be a real vector with $n \geq 1$ components and let $f: \mathbb{R}^{n} \rightarrow \mathbb{R}$ be a smooth function. Then, the unconstrained optimization problem is defined as
\begin{equation}
\min_{\vec{x}\in\mathbb{R}^{n}} f(\vec{x}).
\end{equation}

\subsubsection{Optimality conditions}
For a one-dimensional function $f(x)$ defined and differentiable
over an interval $(a,b)$, the necessary condition for a point $x^* \in (\vec{a},\vec{b})$ to be a local maximum or minimum is that $f'(x^*) = 0$. This is also known as Fermat’s theorem. The multidimensional extension of this condition states that the gradient must be zero at local optimum point, namely
\begin{equation} \label{eq:2}
\nabla f(\vec{x}^*) = 0.
\end{equation}
Eq. \ref{eq:2} is referred as a first order optimality condition, as it expresses in terms of the first order derivatives.

\subsection{Constrained optimization problem}   

\subsubsection{Statement}
The constrained optimization problem is formulated as
\begin{equation}
\min_{\vec{x}\in\mathbb{R}^{n}} f(\vec{x})
\end{equation}
subject to
\begin{equation} \label{eq:4}
c_i\left(\vec{x}\right)\leq 0, \quad i = 1,2,...,w,
\end{equation}
\begin{equation} \label{eq:5}
c_j(\vec{x})=0, \quad j = 1,2,...,l,
\end{equation}
With $c(\vec{x})$ a smooth real-valued function on a subset of $\mathbb{R}^n$. Notably, $c_i(\vec{x})$ and $c_j(\vec{x})$ represent the sets of equality constraints and inequality constraints, respectively. The feasible set is identified as the set of points $\vec{x}$ that satisfy just the constraints (Eqs. \ref{eq:4}, \ref{eq:5}). It must be pointed out that some of the solvers considered in this study are only able to consider equality constraints. In these instances, we will introduce a set of slack variables $s_i$, and convert Eq. \ref{eq:5} into the following set of equality constraints 
\begin{equation}
c_i\left(\vec{x}\right) + s_i^2 = 0, \quad i = 1,2,...,w.
\end{equation}
Such necessary expedient will obviously induce more computational burden on the particular solvers affected by this constraint-type limitation.
%If either the objective n the optimization of nonlinear problems, the constraints functions are given. Most real systems models are non-linear and are much more difficult to optimize compared to other systems model \citep{Chinneck2015}. There are several explanations why this kind of problem is difficult to optimize than the others. Some of them are stated here as follows \citep{Chinneck2015}. 

\subsubsection{Optimality conditions}
The measure of first-order optimality for constrained problems derives from the Karush-Kuhn-Tucker (KKT) conditions \citep{boydlieven}. These necessary conditions are defined as follow. 
Let the objective function $f$ and the constraint functions $g_i$ and $h_j$ be continuously differentiable functions at $\vec{x}^*\in\mathbb{R}^{n}$. If $\vec{x}^*$ is a local optimum and the optimization problem satisfies some regularity conditions \citep{Wright2006}, then there exist the two constants $\mu_i\ (i = 1,\ldots,w)$ and $\lambda_j\ (j = 1,\ldots,\ell)$, called KKT multipliers, such that the following four groups of conditions hold:
\begin{itemize}
    \item \textbf{Stationarity}:
    % \begin{equation}
    % \text{For maximizing} \quad f(\vec{x}): \nabla f(\vec{x}^*) - \sum_{i=1}^m \mu_i \nabla g_i(\vec{x}^*)-\sum_{j=1}^\ell \lambda_j \nabla h_j(\vec{x}^*) = \mathbf 0.
    % \end{equation}
    % \begin{equation}
    %  \text{For minimizing} \quad f(\vec{x}): \nabla f(\vec{x}^*) + \sum_{i=1}^m \mu_i \nabla g_i(\vec{x}^*) + \sum_{j=1}^\ell \lambda_j \nabla h_j(\vec{x}^*) = \mathbf 0.
    % \end{equation}
    \begin{equation}
     f(\vec{x}): \nabla f(\vec{x}^*) + \sum_{i=1}^m \mu_i \nabla g_i(\vec{x}^*) + \sum_{j=1}^\ell \lambda_j \nabla h_j(\vec{x}^*) = \mathbf 0.
    \end{equation}
    \item \textbf{Primal feasibility}:
    \begin{equation}
    g_i(\vec{x}^*) \le 0, \text{ for } i = 1, \ldots, w.
    \end{equation}
    \begin{equation}
    h_j(\vec{x}^*) = 0, \text{ for } j = 1, \ldots, \ell.
    \end{equation}
    \item \textbf{Dual feasibility}:
    \begin{equation}
    \mu_i \ge 0, \text{ for } i = 1, \ldots, w.
    \end{equation}
    \item \textbf{Complementary slackness}:
    \begin{equation}
    \sum_{i=1}^m \mu_i g_i (\vec{x}^*) = 0.
    \end{equation}
\end{itemize}

\section{Selection of NLP solvers and algorithms}
\label{section:s3}
The selection of the NLP solvers considered in this work is based on the following aspects. First of all, we are only considering algorithms that can be implemented in MATLAB. Secondly, we have included solvers that are either free source or, for commercial software, have a free trial version. The remaining part of this section briefly describes the 23 solvers included in our analysis and the most direct source to each algorithm. 

\subsection{APSO}
The Accelerated Particle Swarm Optimization (APSO) is an algorithm developed by Yang at Cambridge University in 2007 and it based on swarm-intelligent search of the optimum  \citep{apso}. APSO is an evolution of the standard particle swarm optimization, and developed to accelerate the convergence of the standard version of the algorithm. The standard PSO is characterized by two elements, the swarm, that is the population, and the members of the population, called particles. The search is based on a randomly initialized population that moves in randomly chosen directions. In particular, each particle moves through the searching space, remembers the best earlier positions, velocity, and accelerations of itself and its neighbors. This information is shared among the particles while they dynamically adjust their own position, velocity and acceleration derived from the best position of all particles. The next step starts when all particles have been shifted. Finally, all particles aim to find the global best among all the current best solutions till the objective function no longer improves or after a certain number of iterations \citep{apso}. The standard PSO uses both the current global best and the individual best, whereas the simplified version APSO is able to accelerate the convergence of the algorithm buy using the global best only. Due to the nature of the algorithm, only constrained nonlinear programming problems can be solved. The MATLAB version of the APSO algorithm is provided in \citep{apso}.

\subsection{BARON}
The Branch and Reduced Optimization Navigator (BARON) is a commercial global optimization software that solves both NLP and mixed-integer nonlinear programs (MINLP). BARON uses deterministic global optimization algorithms of the branch and bound search type which, by applying general assumptions, solve the global optimization problem.  It comes with embedded linear programming (LP) and NLP solvers, such as CLP/CBC, IPOPT, FilterSD and FilterSQP. By default, BARON selects the NLP solver and may switch between different NLP solvers during the search according to problem characteristics and solver performance. To refer to the default option, the name BARON (auto) is chosen. Unlike many other NLP algorithms, BARON doesn’t explicitly require the user to provide an initial guess of the solution but leaves this as an option. If a user doesn't provide the initial guess, then the software shrewdly initializes the variables. In this paper, we use the demo version of the software in conjunction with the MATLAB interface which can be retrieved in \citep{BaronDown}. Must be noted that the free demo version is characterized by some limitations, namely, it can only handle problems with up to ten variables, ten constraints, and it doesn't support trigonometric functions. Details and documentations about BARON software are provided in \citep{Taw2003,Baron}.

\subsubsection{CLP/CBC} 
The Computational Optimization Infrastructure for Operations Research (COIN-OR) Branch and Cut (CBC) is an open-source mixed integer nonlinear programming solver based on the COIN-OR LP solver (CLP) and the COIN-OR Cut generator library (Cgl). The code has been written primarily by John J. Forrest \citep{coinor}.

\subsubsection{IPOPT} 
COIN-OR Interior Point Optimizer (IPOPT) is an open-source solver for large-scale NLP and it has been mainly developed by Andreas Wächter \citep{ipopt}. IPOPT implements an interior point line search filter method for nonlinear programming models. The problem function are not required to be convex but should be twice continuously differentiable. Mathematical details of the algorithm and documentation can be found in \citep{ipopt2}.

\subsubsection{FilterSD} 
FilterSD is a package of Fortran 77 subroutines for solving nonlinear programming problems and linearly constrained problems in continuous optimization. The NLP solver filterSD aims to find a solution of the NLP problem, where the objective function and the constraint function are continuously differentiable at points that satisfy the bounds on $\vec{x}$. The code has been developed to avoid the use of second derivatives, and to prevent storing an approximate reduced Hessian matrix by using a new limited memory spectral gradient approach based on Ritz values. The basic approach is that of Robinson’s method, globalised by using a filter and trust region \citep{filterSD}. 

\subsubsection{FilterSQP} 
FilterSQP is a Sequential Quadratic Programming solver suitable for solving large, sparse or dense linear, quadratic and nonlinear programming problems. The method implements a trust region algorithm with a filter to promote global convergence. The filter accepts a trial point whenever the objective or the constraint violation is improved compared to all previous iterations. The size of the trust region is reduced if the step is rejected, and increased if it is accepted \citep{filterSQP}.

\subsection{FMINCON}
\label{subsection:fmincon}
FMINCON is a MATLAB optimization toolbox used to solve constrained nonlinear programming problems. FMINCON provides the user the option to select amongst five different algorithms to solve nonlinear problems: Active-set, Interior-point, Sequential Quadratic Programming, Sequential Quadratic Programming legacy, and Trust region reflective. Four, out of the five, algorithms are implemented in our analysis as one of them, the Trust Region Reflective algorithm, does not accept the type of constraint considered in our benchmark cases. 

\subsubsection{Active-set}
The Active-set, unlike the Interior point (mentioned next), doesn’t use a barrier term to ensure that the inequality constraints are met, but it solves the optimal equation by understanding the true active-set. A general active-set algorithm for convex quadratic programming can be found in \citep{Wright2006}.

\subsubsection{Interior point}
This method, also known as barrier method, is one type of nonlinear problem-solving algorithms that achieves the determination of optimum values by iteratively approaching the optimal solution from the interior of the feasible set \citep{Wright2006}. Since interior point algorithm depends on a feasible set, the following requirements must be met for the method to be used:
\begin{itemize}
\item the set of feasible interior point should not be empty;
\item all the iterations should occur in the interior of this feasible set.
\end{itemize}

\subsubsection{Sequential Quadratic Programming}
The basic idea behind the Sequential Quadratic Programming (SQP) is to find a minimizer for a subproblem, which is generated as an approximate model of the optimization problem at the current iteration point. This will be then used to define a new iteration point, which is in turn used to define another minimizer, and the process is iterated. SQP is similar to the active set, but some of the differences are listed as follows:
\begin{itemize}
\item strict feasibility with respect to bounds;
\item robustness to non-double results;
\item refactored linear algebra routines;
\item reformulated feasibility routines. 
\end{itemize}
A general line-search algorithm framework for SQP can be found in \citep{Wright2006}.

\subsubsection{Sequential Quadratic Programming legacy}
Sequential Quadratic Programming legacy (SQP-legacy) is similar to SQP, with the difference of using a larger memory and, therefore, it is slower to determining the problem solution \citep{Wright2006}. 

\subsubsection{Trust region reflective}
The Trust region reflective algorithm solves a NLP by defining a region that is assumed to represent the objective function as accurately as possible. From the selected trust region, a step is taken and is used as a minimizer. If that specific step doesn’t generate a solution, a different region with a reduced size will be selected. Then a new step is executed and considered as new minimizer for the region, and the process is iterated \citep{Wright2006}. Trust region reflective algorithm used by FMINCON requires only bounds or linear equality constraints. Due to this setback, this algorithm isn’t included in the analysis. 

\subsection{FMINUNC}
FMINUNC is another MATLAB optimization toolbox used to solve unconstrained nonlinear programming problems \citep{Fminunc}. In this case, FMINUNC gives the user the option of choosing between two different algorithms to solve nonlinear minimization problems: Quasi-Newton, and Trust region.

\subsubsection{Quasi-Newton}
The Quasi-Newton methods build up curvature information at each iteration to formulate a quadratic model problem, with the optimal solution occurring when the sationarity conditions are satisfied. Newton-type methods, as opposed to quasi-Newton methods, calculate the Hessian matrix directly, and proceed in a direction of descent to locate the minimum after a number of iterations, numerically involving a large amount of computation. On the contrary, quasi-Newton methods adopt the observed behavior of the objective function and its gradient to build up curvature information to make an approximation to the Hessian matrix using an appropriate updating technique \citep{Fminunc-qn}. In particular, the quasi-Newton algorithm uses the formula of Broyden, Fletcher, Goldfarb, and Shanno (BFGS) to implement a cubic line search procedure, and for updating the approximation of the Hessian matrix \citep{Fminunc}. 

\subsubsection{Trust region}
The trust region algorithm is a subspace trust-region method, based on the interior-reflective Newton method. Each iteration involves the approximate solution of a large linear system using the method of preconditioned conjugate gradients (PCG) \citep{Fminunc}. In a minimization context, the Hessian matrix can be assumed symmetric, but it is guaranteed to be positive definite only in the neighborhood of a strong minimizer. Algorithm PCG exits when it encounters a direction of negative or zero curvature. The PCG output direction is either a direction of negative curvature or an approximate solution to the Newton system, anyways helping to define the two-dimensional subspace used in the trust-region approach \citep{Fminunc-qn}. 

\subsection{GCMMA}
GCMMA, the Globally Convergent Method of Moving Asymptotes, is a modified version of the MMA that evaluates the global optimum value. However, unlike MMA, GMMA consists of the so called inner and outer iterations. The GCMMA follows the same steps as the MMA except for small changes. In GMMA, an approximate subproblem is created for the first outer iteration by replacing the function with convex functions. These subproblems are then solved in order to find the next iteration points otherwise the inner iteration kicks off. For the first inner iteration, a new subproblem will be generated. Then the next iteration points will be calculated by solving these convex subproblems. The algorithm then moves to the next iteration \citep{gcmma}. The GCMMA algorithm is fully described in \citep{mma2}, and the MATLAB code is freely available at \citep{mma}.

\subsection{KNITRO}
ARTELYS KNITRO is a commercially available nonlinear optimization software package developed by Zienna Optimization since 2001 \citep{Knitro}. KNITRO, short for Nonlinear Interior point Trust Region Optimization, is a software package for finding local solutions of both continuous optimization problems, with or without constraints, and discrete optimization problems with integer or binary variables. The KNITRO package provides efficient and robust solution of small or large problems, for both continuous and discrete problems, derivative-free options. It supports the most popular operating systems and several modeling language and programmatic interfaces \citep{Knitro2}. Multiple versions of the software are available to download at \citep{Knitro}. In this work, the software free trial license is used, in conjunction with the MATLAB interface. Several algorithms are included in the software, such as Interior point, Active-set, and Sequential Quadratic Programming. The description of these algorithms can be found in Section \ref{subsection:fmincon}.

\subsection{MIDACO}
The Mixed Integer Distributed Ant Colony Optimization (MIDACO) is a global optimization solver that combines an extended evolutionary probabilistic technique, called the Ant Colony Optimization algorithm, with the Oracle Penalty method for constrained handling \citep{Midaco}. Ant Colony optimization is modelled using the behavior of ants to find the quickest path between their colony and the source food.
Like the majority of evolutionary optimization algorithms, MIDACO considers the objective and constraint functions as black-box functions. MIDACO was created in collaboration of European Space Agency and EADS Astrium to solve constrained mixed-integer non-linear (MINLP) space applications \citep{midaco2}. We use the trial version of MIDACO, in conjunction with the MATLAB interface. The trial version has a limitation, namely, it doesn't support more than four variables per problem. The solver can be downloaded from \citep{Midaco}.

\subsection{MMA}
The Method of Moving Asymptotes (MMA) solves nonlinear problem function by generating an approximate subproblem. These convex functions used as subproblems are chosen using gradient information at the current iteration points, and also at parameters that are updated at each iteration stage, called the moving asymptotes. The subproblem is solved at the current iteration point, and the solution is used as the next iteration point. Similarly, a new subproblem is generated at this new iteration point, which again is solved to create the next iteration point \citep{mma3}. The MMA algorithm is fully described in \citep{mma2}, and the MATLAB code is freely available at \citep{mma}.

\subsection{MQA}
\label{subsection:s3.9}
The Modified Quasilinearization Algorithm (MQA) is the modified version of the Standard Quasilinearization Algorithm (SQA) \citep{eloe2019quasilinearization,yeo1974quasilinearization} described below. These quasilinearization algorithms base their solution search on the linear approximation of the NLP, namely, on the Hessian matrix and gradient of the objective and constraint functions. Ultimately, the goal is the progressive reduction of the performance index. For unconstrained NLP problems, the performance index is defined as $\Tilde{Q}=f^T_x f_x$, where $f_x$ is the gradient of the objective function. On the other hand, for constrained NLP problems the performance index is defined as $\Tilde{R}=\Tilde{P}+\Tilde{Q}$, which comprises both the feasibility index $\Tilde{P}= h^T h$, and optimality index $\Tilde{Q}=F^T_x F_x$, with $F = f + \lambda^T h$, where $f$ is the objective function, $h$ is the constraint function, and $\lambda$ is the vector of Lagrange multipliers associate to the constraint function. Convergence to the desired solution is achieved when the performance index $\Tilde{Q} \leq \varepsilon_1$ or $\Tilde{R} \leq \varepsilon_2$, with $\varepsilon_1$ and $\varepsilon_2$ small preselected positive constants, for the unconstrained and constrained case respectively \citep{miele1971modified,miele1974modified}. Unlike SQA, characterized by a unitary step size, MQA reduces progressively the step size $0 < \alpha < 1$ to enforce an improvement in optimality. In turn, the main advantage of the MQA, over the SQA, is its descent property: if the stepsize $\alpha$ is sufficiently small, the reduction in the performance index is indeed guaranteed. It must be pointed out that the MQA for NLP problems can only treat equality constraints. Therefore, in our implementation, all the inequality constraints are converted into equality constraints by introducing the slack variables. We have implemented the algorithm on MATLAB in order to solve both unconstrained and constrained NLP problems, based on \citep{miele1971modified,miele1974modified}.

\subsection{PENLAB}
PENLAB is a free open source software package implemented in MATLAB for nonlinear optimization, linear and nonlinear semidefinite optimization and any combination of these. It derives from PENNON, the original implementation of the algorithm which is not open source \citep{penlab}. Originally, PENNON was an implementation of the PBM method developed by Ben-Tal and Zibulevsky for problems of structural optimization, that has grown into a stand alone program for solving general problems \citep{pennon}. It is based on a generalized Augmented Lagrangian method pioneered by R. Polyak \citep{pennon2}. PENLAB can be freely downloaded from \citep{penlabcode}.

\subsection{SGRA}
\label{subsection:s3.11}
The Sequential Gradient-Restoration Algorithm (SGRA) is a first order nonlinear programming solver developed by Angelo Miele and his research group in 1969 \citep{Coker1984,Miele1969}. It is based on a cyclical scheme whereby, first, the constraints are satisfied to a prescribed accuracy (restoration phase); then, using a first-order gradient method, a step is taken toward the optimal direction to improve the performance index (gradient phase). The performance index is defined as $\Tilde{R}=\Tilde{P}+\Tilde{Q}$, which includes both the feasibility index $\Tilde{P}= h^T h$, and optimality index $\Tilde{Q}=F^T_x F_x$, with $F = f + \lambda^T h$, where $f$ is the objective function, $h$ is the constraint function, and $\lambda$ is the vector of Lagrange multipliers associate to the constraint function. Convergence is achieved when the constraint error, and the optimality condition error are $\Tilde{P} \leq \varepsilon_1, \quad \Tilde{Q} \leq \varepsilon_2$, respectively, with $\varepsilon_1$, $\varepsilon_2$ small preselected positive constants. It must be pointed out that the SGRA for NLP problems can only treat equality constraints. Therefore, in our implementation, all the inequality constraints are converted into equality constraints by introducing the slack variables. We have programmed the algorithm on MATLAB in order to solve both unconstrained and constrained NLP problems, based on \citep{Miele1969}. The SGRA version used to solve unconstrained NLP problems differs from the original formulation by the omission of the restoration phase in the iterative process. 

\subsection{SNOPT}
The Sparse Nonlinear OPTimizer (SNOPT) is a commercial software package for solving large-scale optimization problems, linear and nonlinear programs. It minimizes a linear or nonlinear function subject to bounds on the variables and sparse linear or nonlinear constraints. SNOPT implements a sequential quadratic programming method for solving constrained optimization problems with functions and gradients that are expensive to evaluate, and with smooth nonlinear functions in the objective and constraints \citep{snopt}. SNOPT is implemented in Fortran 77 and distributed as source code. In this paper, we use the free trial version of the software in conjunction with the MATLAB interface, that can be retrieved at \citep{snopt2}.

\subsection{SOLNP}
SOLNP is originally implemented in MATLAB to solve general nonlinear programming problems, characterized by nonlinear smooth functions in the objective and constraints \citep{solnp}. Inequality constraints are converted into equality constraints by means of slack variables. The major iteration of SOLNP solves a linearly constrained optimization problem with an augmented Lagrangian objective function. Within the major iteration, as first step it is checked if the solution is feasible for the linear equality constraints of the objective function; if it is not, an interior linear programming procedure is called to find an interior feasible (or near-feasible) solution. Successively, a sequential quadratic programming (QP) solves the linearly constrained problem. If the QP solution is both feasible and optimal, the algorithm stops, otherwise it solves another QP problem as minor iteration. Both major and minor processes repeat until the optimal solution is found or the user-specified maximum number of iterations is reached \citep{solnp}. The SOLNP module in MATLAB can be freely downloaded from \citep{solnpcode}.

\subsection{SQA}
\label{subsection:s3.14}
The Standard Quasilinearization Algorithm (SQA) is the standard version of the QA, and it uses QA techniques for solving nonlinear problems by generating a sequence of linear problems solutions \citep{eloe2019quasilinearization,yeo1974quasilinearization}. SQA differs from MQA for the value associated to the scaling factor $\alpha$. As mentioned before, the SQA can only treat equality constraints. Therefore, in our implementation, all the inequality constraints are converted into equality constraints by introducing the slack variables. 

\section{Convergence metrics and solvers implementation}
\label{section:s4}
In this section we provide the description of the convergence metrics considered in our analysis and narrate the key implementation steps for each solver.

\subsection{Convergence metrics}
The main goal of this paper is to characterize the convergence performance, in terms of speed and accuracy, of the different solvers under analysis. We have selected a number of benchmark NLPs and compared the numerical solutions returned by each solver with the true analytical solution. Moreover, considering that the choice of the initial guesses critically affects the convergence process, we want to also assess the capability to converge to the true optimum, rather than converging into to local minima or not converging at all. With this in mind, we define as converging robustness the quality of a solver to achieve the solution when the search process is initiated from a broad set of initial guesses randomly chosen within the search domain. Finally, to have an accurate assessment of the convergence speed, we will require the solver to repeat the same search several times and average out the total CPU time. As result, given $N$ benchmark test functions, $M$ solvers/algorithms, $K$ randomly generated initial guesses, and $Z$ repeated identical search runs, a total of $N \times M \times K \times Z$ runs have been executed. 

The following performance metrics are in order:

\begin{itemize}
    \item Mean error [\%]:
    \begin{equation}
      \Bar{E}_{m} = \frac{1}{N} \sum_{n=1}^{N} \Bar{E}_{n}, \quad \Bar{E}_{n} = \frac{1}{K} \sum_{k=1}^{K} E_{k}, \quad E_{k}  = 100 \frac{\abs{ f(\vec{x}) - f(\vec{x}^*) }}{\max(\abs{f(x^*)},0.001)}
    \end{equation}
    with $f(\vec{x})$ the benchmark test function evaluated at the numerical solution $\vec{x}$ provided by the solver, $f(\vec{x^*})$ the benchmark test function evaluated at the optimal solution $f(\vec{x^*})$, $E_{k}$ the error associated to the run from the $k$-th randomly generated initial guess, $\Bar{E}_{n}$ the mean error associated to the $n$-th benchmark test function, and $\Bar{E}_{m}$ the mean error delivered by the $m$-th solver. The biunivocal choice of the denominator of $E_{k}$ is based on the fact that some benchmark test functions at the optimal solution have zero value; in this case, a value of $0.001$ is chosen instead as reference value.
    
    \item Mean variance [\%]:
    \begin{equation}
      \Bar{\sigma}_{m} = \frac{1}{N} \sum_{n=1}^{N} \sigma_{n}, \quad \sigma_{n}  =  \frac{1}{K-1} \sum_{k=1}^{K} \abs{E_{k} - \Bar{E}_{n}}^2
    \end{equation}
    where $\sigma_{n}$ is the variance correspondent to the $n$-th benchmark test function, and $\Bar{\sigma}_{m}$ the mean variance delivered by the $m$-th solver. 
    
    \item Mean convergence rate [\%]:
    \begin{equation}
        \Bar{\gamma}_{m} = \frac{1}{N} \sum_{n=1}^{N} \gamma_{n}, \quad \gamma_{n} = 100 \frac{K-K_{conv}}{K}
    \end{equation}
    with $K_{conv}$ the number of runs (from a pool of $K$ distinct initial guesses) which successfully reach convergence for the $n$-th function, $\gamma_{n}$ the convergence rate for the $n$-th function, and $\Bar{\gamma}_{m}$ the mean convergence rate delivered by the $m$-th solver. The convergence rate is computed considering succesfull a run that satisfies the converging threshold conditions $E_k \leq E_{max}=5\%$, and $CPU_k \leq CPU_{max} = 10$ s, with $CPU_k$ is the CPU time required to the run starting from the $k$-th initial guess. 
    
    \item Mean CPU time [s]:
    \begin{equation}
      \overline{CPU}_{m} = \frac{1}{N} \sum_{n=1}^{N}  \overline{CPU}_{n},
    \end{equation}
    \begin{equation}  
      \overline{CPU}_{n} = \frac{1}{Z} \sum_{z=1}^{Z} \overline{CPU}_{z}, \quad  \overline{CPU}_{z} = \frac{1}{K} \sum_{k=1}^{K} CPU_k
    \end{equation}
    where  $\overline{CPU}_{z}$ is the mean CPU time per $z$-th repetition, $\overline{CPU}_{n}$ is the mean CPU time related to the $n$-th benchmark test function, and $\overline{CPU}_{m}$ is the mean CPU time delivered by the $m$-th solver.
\end{itemize}

\subsection{Solvers implementation}
\label{subsection:s4.2}
In this paper we analyze the convergence performances of the different solvers in terms of robustness, accuracy, and convergence speed. Considering that the user might decide to tune the convergence parameters to favor one of these metrics, we have decided to perform the comparison for three separate implementation scenarios: plug and play (P\&P), high accuracy (HA), and quick solution (QS). The plug and play settings, as the name suggests, are the "out-of-the-box" settings of each solver. The high accuracy settings are based on more stringent tolerances and/or on a higher number of maximum iterations with respect to the plug and play settings. This tuning aims to achieve a more precise solution. Finally, the quick solution settings are characterized by more relaxed convergence tolerances, and a lower number of maximum iterations with respect to the plug and play settings. In this scenario the algorithms should reach a less accurate solution but in a shorter time. In general, the objective function, its gradient, the initial conditions, the constraint function (for constrained problems only), and the solver options are elements which are inputted to each solver. The objective function gradient is not necessary for APSO, BARON, MIDACO, and SOLNP, but it is optional for FMINCON/FMINUNC and KNITRO. For GCMMA/MMA, SGRA, and SNOPT, the gradient of both the objective and constraint functions is necessary. MQA/SQA and PENLAB, in addition to these inputs, require the Hessian of the objective function. In the following subsection, details on each solver, and on the three different solver settings per each solver are described. It must be noted that, in most cases, the settings' names here reported are the same of the solver's options names used in the code implementation. In this way, the reader can have a better understanding of which solver's parameter has been tuned. 

\subsubsection{APSO}
The three settings considered in the analysis are reported in Table \ref{tab:apso_sett}, where \textit{no. particles} is the number of particles, \textit{no. iterations} is the total number of iterations, and $\gamma$ is a control parameter that multiplies $\alpha$, one of the two learning parameters or acceleration constants, $\alpha$ and $\beta$, the random amplitude of roaming particles and the speed of convergence, respectively. APSO does also require the number of problem variables, \textit{no. vars}, to be defined but this parameter is, obviously, invariant for the three settings.

\begin{table}[ht]
	\caption{APSO settings.}
	\label{tab:apso_sett}
% 	\customsize
	\small
% 	\footnotesize
	\centering
	\begin{tabular}{lccc}
		\toprule
		Settings & P\&P & HA & QS  \\ 
		\midrule
        \textit{no. particles} &  15 & 50 & 10 \\
        \textit{no. iterations} & 300 & 500 & 100 \\
        $\gamma$ & 0.9 & 0.95 & 0.95 \\
		\bottomrule    
	\end{tabular}
\end{table}

\subsubsection{BARON}
The three settings considered in the analysis are reported in Table \ref{tab:baron_sett}, with $EpsA$ the absolute termination tolerance, $EpsR$ the relative termination tolerance, and $AbsConFeasTol$ the absolute constraint feasibility tolerance. Due to the limitations of the trial version of the solver, trigonometric functions and problems with more than ten variables are not supported by the solver; for this reason, the following test functions are excluded in the analysis: A.2, A.3, A.4, A.5, A.7, A.11, A.13, A.14, A.16, A.17, A.18, A.22, A.24, A.26 for unconstrained problems, and B.1, B.2, B.5, B.8, B.20 for constrained problems.

\begin{table}[ht]
	\caption{BARON settings.}
	\label{tab:baron_sett}
% 	\customsize
	\small
% 	\footnotesize
	\centering
	\begin{tabular}{lccc}
		\toprule
		Settings & P\&P & HA & QS  \\ 
		\midrule
        $EpsA$ & 1e-6 & 1e-10 & 1e-3 \\
        $EpsR$ & 1e-4 & 1e-10 & 1e-3 \\
        $AbsConFeasTol$ & 1e-5 & 1e-10 & 1e-3 \\
		\bottomrule    
	\end{tabular}
\end{table}

\subsubsection{FMINCON/FMINUNC}
The three settings considered in the analysis are reported in Table \ref{tab:fmin_sett}, with $StepTolerance$ the lower bound on the size of a step, $ConstraintTolerance$ the upper bound on the magnitude of any constraint functions, $FunctionTolerance$ the lower bound on the change in the value of the objective function during a step, and $OptimalityTolerance$ the tolerance for the first-order optimality measure. 

\begin{table}[ht]
	\caption{FMINCON/FMINUNC settings.}
	\label{tab:fmin_sett}
% 	\customsize
	\small
% 	\footnotesize
	\centering
	\begin{tabular}{lccc}
		\toprule
		Settings & P\&P & HA & QS  \\ 
		\midrule
		FMINCON \\
        $StepTolerance$ & 1e-10 & 1e-10 & 1e-6 \\
        $ConstraintTolerance$ & 1e-6 & 1e-10 & 1e-3 \\
        $FunctionTolerance$ & 1e-6 & 1e-10 & 1e-3 \\
        $OptimalityTolerance$ & 1e-6 & 1e-10 & 1e-3 \\
        \midrule
		FMINUNC (quasi-newton) \\
        $StepTolerance$ & 1e-6 & 1e-12 & 1e-6 \\
        $FunctionTolerance$ & 1e-6 & 1e-12 & 1e-3 \\
        $OptimalityTolerance$ & 1e-6 & 1e-12 & 1e-3 \\
        \midrule
        FMINUNC (trust-region) \\
        $StepTolerance$ & 1e-6 & 1e-12 & 1e-6 \\
        $FunctionTolerance$ & 1e-6 & 1e-12 & 1e-3 \\
        $OptimalityTolerance$ & 1e-6 & 1e-6 & 1e-3 \\
		\bottomrule    
	\end{tabular}
\end{table}

\subsubsection{GCMMA/MMA}
The three settings considered in the analysis are reported in Table \ref{tab:mma_sett}, where $epsimin$ is a prescribed small positive tolerance that terminates the algorithm, whereas $maxoutit$ is the maximum number of iterations for MMA, and the maximum number of outer iterations for GCMMA.

\begin{table}[ht]
	\caption{GCMMA/MMA settings.}
	\label{tab:mma_sett}
% 	\customsize
	\small
% 	\footnotesize
	\centering
	\begin{tabular}{lccc}
		\toprule
		Settings & P\&P & HA & QS  \\ 
		\midrule
        $epsimin$ & 1e-7 & 1e-10 & 1e-3 \\
        $maxoutit$ & 80 & 150 & 30\\
		\bottomrule    
	\end{tabular}
\end{table}

\subsubsection{KNITRO}
The three settings considered in the analysis are reported in Table \ref{tab:knitro_sett}, where $MaxIter$ is the maximum number of iterations before termination, $TolX$ is a tolerance that terminates the optimization process if the relative change of the solution point estimate is less than that value, $TolFun$ specifies the final relative stopping tolerance for the KKT (optimality) error, and $TolCon$ specifies the final relative stopping tolerance for the feasibility error.

\begin{table}[ht]
	\caption{KNITRO settings.}
	\label{tab:knitro_sett}
% 	\customsize
	\small
% 	\footnotesize
	\centering
	\begin{tabular}{lccc}
		\toprule
		Settings & P\&P & HA & QS  \\ 
		\midrule
        $MaxIter$ & 1000 & 10000 & 100\\
        $TolX$ & 1e-6 & 1e-10 & 1e-3\\
        $TolFun$ & 1e-6 & 1e-10 & 1e-3\\
        $TolCon$ & 1e-6 & 1e-10 & 1e-3\\
		\bottomrule    
	\end{tabular}
\end{table}

\subsubsection{MIDACO}
The three settings considered in the analysis are reported in Table \ref{tab:midaco_sett}, where $maxeval$ is the maximum number of function evaluation. It is a distinctive feature of MIDACO that allows the solver to stop exactly after that number of function evaluation. Due to the limitations of the trial version of the solver, test functions with more than four variables are not supported by the solver; for this reason, the following test functions are excluded in the analysis: A.7, A.10, A.11, A.13, A.14, A.15, A.16, A.18, A.22, A.24 for unconstrained problems, and B.1, B.2, B.3, B.4, B.7, B.9, B.10, B.12, B.18, B.19, B.20, B.21, B.22, B.26, B.29, B.30 for constrained problems.

\begin{table}[ht]
	\caption{MIDACO settings.}
	\label{tab:midaco_sett}
% 	\customsize
	\small
% 	\footnotesize
	\centering
	\begin{tabular}{lccc}
		\toprule
		Settings & P\&P & HA & QS  \\ 
		\midrule
        $maxeval$ & 50000 & 150000 & 10000 \\
		\bottomrule    
	\end{tabular}
\end{table}

\subsubsection{MQA}
The three settings considered in the analysis are reported in Table \ref{tab:mqa_sett}, with $\varepsilon_1$ and $\varepsilon_2$  the prescribed small positive tolerances that allow the solver to stop, when the inequality $\Tilde{Q} \leq \varepsilon_1$ or $\Tilde{R} \leq \varepsilon_2$ is met. As mentioned in Section \ref{subsection:s3.9}, MQA for NLP problems can only treat equality constraints, namely all the inequality constraints are converted into equality constraints by introducing the slack variables. In this study, for all the three settings considered in the analysis, a value of $1$ is chosen as initial guess for all the slack variables.

\begin{table}[ht]
	\caption{MQA settings.}
	\label{tab:mqa_sett}
% 	\customsize
	\small
% 	\footnotesize
	\centering
	\begin{tabular}{lccc}
		\toprule
		Settings & P\&P & HA & QS  \\ 
		\midrule
        $\varepsilon_1$ & 1e-5 & 1e-8 & 1e-2 \\
        $\varepsilon_2$ & 1e-4 & 1e-5 & 1e-3 \\
		\bottomrule    
	\end{tabular}
\end{table}

\subsubsection{PENLAB}
The three settings considered in the analysis are reported in Table \ref{tab:penlab_sett}, where $max\_inner\_iter$ is the maximum number of inner iterations, $max\_outer\_iter$ is the maximum number of outer iterations, $mpenalty\_min$ is the lower bound for penalty parameters, $inner\_stop\_limit$ is the termination tolerance for the inner iterations, $outer\_stop\_limit$ is the termination tolerance for the outer iterations, $kkt\_stop\_limit$ is the termination tolerance KKT optimality conditions, and $unc\_dir\_stop\_limit$ is the stopping tolerance for the unconstrained minimization.

\begin{table}[ht]
	\caption{PENLAB settings.}
	\label{tab:penlab_sett}
% 	\customsize
	\small
% 	\footnotesize
	\centering
	\begin{tabular}{lccc}
		\toprule
		Settings & P\&P & HA & QS  \\ 
		\midrule
        $max\_inner\_iter$ & 100 & 1000 & 25 \\
        $max\_outer\_iter$ & 100 & 1000 & 25 \\
        $mpenalty\_min$ & 1e-6 & 1e-9 & 1e-3 \\
        $inner\_stop\_limit$ & 1e-2 & 1e-9 & 1e-1 \\
        $outer\_stop\_limit$ & 1e-6 & 1e-9 & 1e-3 \\
        $kkt\_stop\_limit$ & 1e-4 & 1e-6 & 1e-2 \\
        $unc\_dir\_stop\_limit$ & 1e-2 & 1e-9 & 1e-1 \\
		\bottomrule    
	\end{tabular}
\end{table}

\subsubsection{SGRA}
The three settings considered in the analysis are reported in Table \ref{tab:sgra_sett}, with $\varepsilon_{1}$ the tolerance related to the constraint error $\Tilde{P}$, and $\varepsilon_{2}$ the tolerance related to the optimality condition error $\Tilde{Q}$. Considering that the SGRA can only treat equality constraints, all the inequality constraints are converted into equality constraints by introducing the slack variables. In this study, for all the three settings considered in the analysis, a value of $1$ is chosen for all the slack variables.

\begin{table}[ht]
	\caption{SGRA settings.}
	\label{tab:sgra_sett}
% 	\customsize
	\small
% 	\footnotesize
	\centering
	\begin{tabular}{lccc}
		\toprule
		Settings & P\&P & HA & QS  \\ 
		\midrule
        $\varepsilon_{1}$ & 1e-9 & 1e-10 & 1e-8 \\
        $\varepsilon_{2}$ & 1e-4 & 1e-6 & 1e-2 \\
		\bottomrule    
	\end{tabular}
\end{table}

\subsubsection{SNOPT}
The three settings considered in the analysis are reported in Table \ref{tab:snopt_sett}, where $major\_iterations\_limit$ is the limit on the number of major iterations in the SQP method, 
$minor\_iterations\_limit$ is the limit on minor iterations in the QP subproblems, $major\_feasibility\_tolerance$ is the tolerance for feasibility of the nonlinear constraints, $major\_optimality\_tolerance$ is the tolerance for the dual variables, and $minor\_feasibility\_tolerance$ is the tolerance for the variables and their bounds.

\begin{table}[ht]
	\caption{SNOPT settings.}
	\label{tab:snopt_sett}
% 	\customsize
	\small
% 	\footnotesize
	\centering
	\begin{tabular}{lccc}
		\toprule
		Settings & P\&P & HA & QS  \\ 
		\midrule
        $major\_iterations\_limit$ & 1000 & 10000 & 100 \\
        $minor\_iterations\_limit$ & 500 & 5000 & 100 \\
        $major\_feasibility\_tolerance$ & 1e-6 & 1e-12 & 1e-3 \\
        $major\_optimality\_tolerance$ & 1e-6 & 1e-12 & 1e-3 \\
        $minor\_feasibility\_tolerance$ & 1e-6 & 1e-12 & 1e-3 \\
		\bottomrule    
	\end{tabular}
\end{table}

\subsubsection{SOLNP}
The three settings considered in the analysis are reported in Table \ref{tab:solnp_sett}, with $\rho$ the penalty parameter in the augmented Lagrangian objective function, $maj$ the maximum number of major iterations, $min$ the maximum number of minor iterations, $\delta$ the perturbation parameter for numerical gradient calculation, and $\epsilon$ the relative tolerance on optimality and feasibility. During the HA scenario implementation, we learned that different convergence settings are required for unconstrained and constrained problems. This peculiarity might be induced by the stringent tolerances adopted in this scenario.

\begin{table}[ht]
	\caption{SOLNP settings. Tuning values for the HA scenario are divided for unconstrained (left-side) and constrained (right-side) problems.}
	\label{tab:solnp_sett}
% 	\customsize
% 	\small
% 	\footnotesize
	\centering
	\begin{tabular}{lccc}
		\toprule
		Settings & P\&P & HA & QS  \\ 
		\midrule
        $\rho$ & 1 & 1 & 1 \\
        $maj$ & 10 & 500$\mid$10 & 10 \\
        $min$ & 10 & 500$\mid$10 & 10 \\
        $\delta$ & 1e-5 & 1e-10$\mid$1e-6 & 1e-3 \\
        $\epsilon$ & 1e-4 & 1e-12$\mid$1e-7 & 1e-3 \\
		\bottomrule    
	\end{tabular}
\end{table}

\subsubsection{SQA}
The three settings considered in the analysis are reported in Table \ref{tab:sqa_sett}, with $\varepsilon_1$ and $\varepsilon_2$  the prescribed small positive tolerances that allow the solver to stop, when the inequality $\Tilde{Q} \leq \varepsilon_1$ or $\Tilde{R} \leq \varepsilon_2$ is met. As mentioned earlier, SQA can only treat equality constraints. To overcome this limitation, the inequality constraints are converted into equality constraints by introducing slack variables. In this study, for all the three settings considered in the analysis, a value of $1$ is chosen for all the slack variables.

\begin{table}[ht]
	\caption{SQA settings.}
	\label{tab:sqa_sett}
% 	\customsize
	\small
% 	\footnotesize
	\centering
	\begin{tabular}{lccc}
		\toprule
		Settings & P\&P & HA & QS  \\ 
		\midrule
        $\varepsilon_1$ & 1e-5 & 1e-8 & 1e-2 \\
        $\varepsilon_2$ & 1e-4 & 1e-5 & 1e-3 \\
		\bottomrule    
	\end{tabular}
\end{table}

\section{Benchmark test functions and results} 
\label{section:s5}
We present a collection of unconstrained and constrained optimization test problems that are used to validate the performance of the various optimization algorithms presented above for the different implementation scenarios. The comparison results are also discussed in depth in this section.

For performance comparison purposes, an equivalent environment and control parameters have been created to run each NLP solver. All outputs tabulated in this paper are calculated using MATLAB software running on a desktop computer with the following specs: Intel(R) Core(TM) i7-6700 CPU 3.40GHz processor, 16.0 GB of RAM, running a 64-bit Windows 10 operating system. To assess the true computational time required by each algorithm to reach convergence, implementation steps that are expected to have an impact on the computer's performance are deactivated during the run for the solution. The internet connection and other unrelated applications are turned off throughout the analysis, ensuring that unnecessary background activities are not accessing computational resources throughout the solvers' performance. A collection of unconstrained and constrained benchmark problems are used to test the solvers based on, but not limited to \citep{testcollect,benchmark,handbook}. Specifically, the benchmark problems include combinations of logarithmic, trigonometric, and exponential terms, non-convex and convex functions, a minimum of two to a maximum of thirty variables, and a maximum of nine constraint functions for the constrained optimization problems. For sake of completeness, all the benchmark test functions are listed in Appendix \ref{sec:appendixA} and \ref{sec:appendixB}. For each of the test function, dimension, domain and search space, objective function, constraints, and minimum solution are listed. As mentioned in Section \ref{subsection:s4.2}, the comparison between each solver is carried out by considering three different settings: plug and play, high accuracy, and quick solution. In this way we want to assess the robustness, accuracy, and convergence speed of every solver.  For each benchmark problem, all solvers use the same set of randomly generated initial guesses.

\subsection{Results for unconstrained optimization problems}
A collection of 30 unconstrained optimization test problems is used to validate the performance of the optimization algorithms. The benchmark test functions for unconstrained global optimization are listed in Appendix \ref{sec:appendixA}. For the purpose of this analysis, given $N = 30$ benchmark test functions, $M = 16$ solvers and algorithms, $K = 50$ randomly generated initial guesses, and $Z = 3$ iterations, a set of $N \times M \times K \times Z$ runs are executed. Tables \ref{tab:prallu1}, \ref{tab:prallu2}, and \ref{tab:prallu3} report the results for the plug and play (P\&P), high accuracy (HA), and quick solution (QS) settings, respectively. From the analysis of the results for the P\&P settings, Table \ref{tab:prallu1}, we observe that BARON (auto) and BARON (ipopt) are able to reach the minimum mean error and variance, the highest convergence rate, but they are not the fastest ones to reach the solution. Moreover, BARON (sd), BARON (sqp), SNOPT, and PENLAB are able to obtain good results in terms of mean error and variance. Overall, PENLAB is also able to reach a convergence rate similar to BARON (auto) and BARON (ipopt), with the advantage of being 10 times faster than them. The worst results in terms of accuracy and convergence rate are obtained by SOLNP and SGRA. For the HA settings, Table \ref{tab:prallu2}, we can observe similar trends. In general, as expected, all the solvers manage to achieve a more accurate solution as they reduce the average error, increase their convergence rate, and increase the average convergence time. MIDACO is now able to reach the highest convergence rate together with all the versions of BARON. Overall PENLAB is the solver which delivers a good trade-off in performance. With respect the P\&P settings, SOLNP significantly improves its convergence rate, whereas SGRA just slightly increase its performances. It is interesting to observe that, KNITRO (interior-point) and KNITRO (sqp), aside improving their convergence rate, increase their mean error and variance increase. Despite our effort, we are not sure how to explain this unexpected behaviour. Regarding the QS settings, Table \ref{tab:prallu3}, generally all the solvers reduce their convergence time and also decrease their convergence rate except for BARON (auto), BARON (ipopt), and BARON (sqp) which remain unaltered. SQA, FMINUNC, SOLNP, and FMINUNC are amongst the fastest to reach the solution but their convergence rate is quite low. In addition, conversely to all the other solvers that experience a smaller CPU time, BARON is not always able to achieve a faster CPU time with respect to the P\&P settings. The same happens to the SGRA, probably due to its intrinsic iterative nature. 

\FloatBarrier
\begin{table}[ht]
	\caption{All unconstrained problems, plug and play (P\&P) settings. Solvers ranked w.r.t. convergence rate.}
	\label{tab:prallu1}
% 	\customsize
 	\small
% 	\footnotesize
	\centering
	\begin{tabular}{clcccc}
		\toprule
		Ranking & Solver & $\Bar{E}$[\%] & $\Bar{\sigma}$[\%] & $\Bar{\gamma}$[\%] & $\overline{CPU}$[$s$]  \\ 
		\midrule
1 & BARON (auto) & 2.559e-07 & 2.171e-31 & 92.3 & 0.3434 \\ 
2 & BARON (ipopt) & 2.559e-07 & 2.162e-31 & 92.3 & 0.2845 \\ 
3 & BARON (sd) & 1.596e-03 & 2.615e-09 & 92.3 & 0.3863 \\ 
4 & BARON (sqp) & 2.536e-07 & 4.231e-20 & 92.3 & 0.3389 \\ 
12 & FMINUNC (quasi-newton) & 8.643e-02 & 1.669e-01 & 52.8 & 0.0045 \\ 
8 & FMINUNC (trust-region) & 7.836e-02 & 2.068e-02 & 68.8 & 0.0153 \\
11 & KNITRO (active-set) & 2.299e-01 & 4.979e-02 & 60.5 & 0.0200 \\
10 & KNITRO (interior-point) & 2.703e-02 & 3.235e-02 & 61.2 & 0.0194 \\
9 & KNITRO (sqp) & 1.121e-02 & 1.274e-02 & 61.7 & 0.0365 \\
6 & MIDACO & 3.130e-01 & 1.626e-01 & 85.2 & 0.3193 \\ 
14 & MQA & 2.031e-01 & 8.748e-02 & 51.7 & 0.1345 \\ 
5 & PENLAB & 1.016e-03 & 5.340e-37 & 88.5 & 0.0125 \\
16 & SGRA & 5.921e-01 & 8.627e-02 & 40.8 & 0.2227 \\
7 & SNOPT & 7.008e-03 & 2.444e-02 & 73.8 & 0.0071 \\
15 & SOLNP & 4.648e-01 & 1.908e-01 & 48.2 & 0.0097 \\ 
13 & SQA & 2.362e-01 & 1.383e-01 & 52.5 & 0.0005 \\ 
		\bottomrule    
	\end{tabular}
\end{table}
\FloatBarrier

\FloatBarrier
\begin{table}[ht]
	\caption{All unconstrained problems, high accuracy (HA) settings. Solvers ranked w.r.t. mean error.}
	\label{tab:prallu2}
% 	\customsize
	\small
% 	\footnotesize
	\centering
	\begin{tabular}{clcccc}
		\toprule
		Ranking & Solver & $\Bar{E}$[\%] & $\Bar{\sigma}$[\%] & $\Bar{\gamma}$[\%] & $\overline{CPU}$[$s$]  \\ 
		\midrule
2 & BARON (auto) & 2.563e-07 & 2.171e-31 & 92.3 & 0.6550 \\ 
3 & BARON (ipopt) & 2.563e-07 & 2.162e-31 & 92.3 & 0.6624 \\ 
4 & BARON (sd) & 1.186e-06 & 3.611e-13 & 92.3 & 0.8777 \\ 
1 & BARON (sqp) & 2.536e-07 & 4.231e-20 & 92.3 & 0.6334 \\ 
9 & FMINUNC (quasi-newton) & 3.526e-03 & 8.852e-04 & 59.2 & 0.0062 \\ 
11 & FMINUNC (trust-region) & 1.860e-02 & 1.423e-02 & 68.8 & 0.0238 \\
13 & KNITRO (active-set) & 7.153e-02 & 1.883e-01 & 67.9 & 0.0273 \\ 
12 & KNITRO (interior-point) & 5.329e-02 & 1.258e-01 & 68.3 & 0.0440 \\ 
14 & KNITRO (sqp) & 9.411e-02 & 1.571e-01 & 69.1 & 0.0731 \\ 
15 & MIDACO & 1.756e-01 & 9.387e-02 & 92.2 & 1.0174 \\ 
8 & MQA & 3.160e-03 & 5.835e-05 & 52.1 & 0.1520 \\
5 & PENLAB & 4.042e-06 & 8.944e-42 & 88.5 & 0.0121 \\
16 & SGRA & 2.709e-01 & 1.335e-01 & 44.9 & 0.2555 \\ 
7 & SNOPT & 1.260e-03 & 1.298e-03 & 74.2 & 0.0099 \\
6 & SOLNP & 9.420e-04 & 7.900e-04 & 69.1 & 0.0095 \\
10 & SQA & 3.984e-03 & 9.000e-05 & 53.2 & 0.0003 \\
		\bottomrule    
	\end{tabular}
\end{table}
\FloatBarrier

\FloatBarrier
\begin{table}[ht]
	\caption{All unconstrained problems, quick solution (QS) settings. Solvers ranked w.r.t. mean CPU time.}
	\label{tab:prallu3}
% 	\customsize
	\small
% 	\footnotesize
	\centering
	\begin{tabular}{clcccc}
		\toprule
		Ranking & Solver & $\Bar{E}$[\%] & $\Bar{\sigma}$[\%] & $\Bar{\gamma}$[\%] & $\overline{CPU}$[$s$]  \\ 
		\midrule
15 & BARON (auto) & 2.556e-07 & 2.171e-31 & 92.3 & 0.3692 \\ 
16 & BARON (ipopt) & 2.556e-07 & 2.162e-31 & 92.3 & 0.3743 \\ 
13 & BARON (sd) & 4.295e-06 & 2.853e-09 & 84.6 & 0.3684 \\ 
14 & BARON (sqp) & 2.536e-07 & 4.231e-20 & 92.3 & 0.3690 \\
2 & FMINUNC (quasi-newton) & 6.076e-01 & 8.157e-01 & 33.8 & 0.0024 \\ 
5 & FMINUNC (trust-region) & 1.924e-01 & 1.997e-01 & 49.3 & 0.0108 \\ 
8 & KNITRO (active-set) & 3.522e-01 & 3.677e-01 & 48.7 & 0.0171 \\ 
7 & KNITRO (interior-point) & 3.231e-01 & 4.066e-01 & 49.1 & 0.0169 \\ 
9 & KNITRO (sqp) & 3.900e-01 & 5.835e-01 & 50.4 & 0.0256 \\
10 & MIDACO & 5.852e-02 & 7.128e-02 & 72.3 & 0.0692 \\ 
11 & MQA & 2.405e-01 & 2.930e-01 & 42.2 & 0.1819 \\ 
6 & PENLAB & 5.452e-05 & 5.623e-39 & 84.6 & 0.0118 \\ 
12 & SGRA & 8.640e-01 & 2.211e-01 & 23.8 & 0.3033 \\ 
3 & SNOPT & 1.581e-01 & 1.367e-01 & 66.4 & 0.0040 \\ 
4 & SOLNP & 5.357e-01 & 3.847e-01 & 41.2 & 0.0093 \\ 
1 & SQA & 1.964e-01 & 1.609e-01 & 43.3 & 0.0002 \\
		\bottomrule    
	\end{tabular}
\end{table}
\FloatBarrier
\FloatBarrier
% \begin{figure}[ht]
% 	\centering
% 	\includegraphics[width=\linewidth]{./figure/mean_error_u.png}
% 	\caption{Unconstrained problems, mean error.}
% 	\label{fig:meanEu}
% \end{figure}

% \begin{figure}[ht]
% 	\centering
% 	\includegraphics[width=\linewidth]{./figure/mean_var_u.png}
% 	\caption{Unconstrained problems, mean variance.}
% 	\label{fig:meanVu}
% \end{figure}

% \begin{figure}[ht]
% 	\centering
% 	\includegraphics[width=\linewidth]{./figure/mean_conv_u.png}
% 	\caption{Unconstrained problems, mean convergence rate.}
% 	\label{fig:meanCVu}
% \end{figure}

% \begin{figure}[ht]
% 	\centering
% 	\includegraphics[width=\linewidth]{./figure/mean_cpu_u.png}
% 	\caption{Unconstrained problems, mean CPU time.}
% 	\label{fig:meanCPUu}
% \end{figure}

\subsection{Results for constrained optimization problems}
A collection of 30 unconstrained optimization test problems is used to validate the performance of the optimization algorithms. The benchmark test functions for constrained global optimization are listed in Appendix \ref{sec:appendixB}. For the purpose of the analysis, given $N = 30$ benchmark test functions, $M = 21$ solvers and algorithms, $K = 50$ randomly generated initial guesses, and $Z = 3$ iterations, a set of $N \times M \times K \times Z$ runs are executed. Tables \ref{tab:prallc1}, \ref{tab:prallc2}, and \ref{tab:prallc3} report the results for the P\&P, HA, and QS settings, respectively.  From the analysis of the results for the P\&P settings, Table \ref{tab:prallc1}, we observe that all the versions of BARON are able to reach almost the highest accuracy and the best convergence rate but they are not the fastest to reach the solution. MIDACO is able to achieve the second best convergence rate, with an average CPU time that is more than 50\% faster than BARON. PENLAB obtains the best mean error and variance but this performance is tempered by a low convergence rate, together with the SGRA, MQA, and SQA which are also quite slow to reach solution. FMINCON (interior-point), KNITRO (interior-point), and SNOPT reach a convergence rate lower than BARON and MIDACO, but they are significantly faster. Regarding the HA settings, Table \ref{tab:prallc2}, similar consideration can be made for BARON, but in this case the CPU time is considerably increasing. MIDACO shows an improvement in the convergence rate, reaching values very similar to BARON. PENLAB still obtains the best mean error and variance, but it has one of the lowest convergence rate, together with the SGRA. In general, most of the solvers increase their convergence rate, and decrease their mean error, except for GCMMA and PENLAB. Regarding the QS settings, Table \ref{tab:prallc3}, generally all the solvers decrease their convergence rate except for BARON and PENLAB. Same considerations about BARON and PENLAB can be done as in the two previous scenarios. MIDACO reports a significant decrease in the convergence rate. The different versions of BARON have similar CPU time with respect to the P\&P settings. FMINCON (interior-point), KNITRO (interior-point), and SNOPT reach a convergence rate lower than BARON, but they are significantly faster. The worst results in terms of convergence rate and CPU time are obtained by MQA and SQA.

\FloatBarrier
\begin{table}[ht]
	\caption{All constrained problems, plug and play (P\&P) settings. Solvers ranked w.r.t. convergence rate.}
	\label{tab:prallc1}
% 	\customsize
	\small
% 	\footnotesize
	\centering
	\begin{tabular}{clcccc}
		\toprule
		Ranking & Solver & $\Bar{E}$[\%] & $\Bar{\sigma}$[\%] & $\Bar{\gamma}$[\%] & $\overline{CPU}$[$s$]  \\ 
		\midrule
17 & APSO & 1.512e+00 & 1.025e+00 & 39.2 & 0.1772 \\ 
1 & BARON (auto) & 2.153e-03 & 2.728e-08 & 92.0 & 0.7016 \\ 
2 & BARON (ipopt) & 2.441e-03 & 2.920e-08 & 92.0 & 0.8052 \\ 
3 & BARON (sd) & 2.162e-03 & 2.566e-08 & 92.0 & 0.7539 \\ 
4 & BARON (sqp) & 2.183e-03 & 1.478e-08 & 92.0 & 0.8188 \\
10 & FMINCON (active-set) & 1.795e-01 & 2.123e-01 & 71.9 & 0.0204 \\ 
6 & FMINCON (interior-point) & 1.985e-01 & 2.413e-01 & 75.9 & 0.0271 \\ 
13 & FMINCON (sqp) & 1.908e-01 & 2.446e-01 & 69.3 & 0.0093 \\ 
11 & FMINCON (sqp-legacy) & 1.893e-01 & 2.429e-01 & 69.4 & 0.0111 \\ 
15 & GCMMA & 4.490e-01 & 3.742e-01 & 45.7 & 0.9681 \\ 
12 & KNITRO (active-set) & 1.908e-01 & 2.759e-01 & 69.4 & 0.0472 \\ 
7 & KNITRO (interior-point) & 1.718e-01 & 1.962e-01 & 74.6 & 0.0303 \\ 
8 & KNITRO (sqp) & 1.788e-01 & 2.027e-01 & 72.9 & 0.1016 \\ 
5 & MIDACO & 4.739e-01 & 2.718e-01 & 81.0 & 0.3331 \\ 
16 & MMA & 7.188e-01 & 5.743e-01 & 44.1 & 0.5856 \\ 
20 & MQA & 5.125e-01 & 3.460e-01 & 20.8 & 3.1559 \\ 
18 & PENLAB & 1.127e-04 & 3.258e-41 & 31.0 & 0.0379 \\ 
19 & SGRA & 6.360e-01 & 7.011e-01 & 30.3 & 0.9815 \\
9 & SNOPT & 1.689e-01 & 2.010e-01 & 72.1 & 0.0040 \\
14 & SOLNP & 3.243e-01 & 3.211e-01 & 48.1 & 0.0095 \\ 
21 & SQA & 3.990e-01 & 5.778e-01 & 20.2 & 3.1822 \\ 
		\bottomrule    
	\end{tabular}
\end{table}
\FloatBarrier

\FloatBarrier
\begin{table}[ht]
	\caption{All constrained problems, high accuracy (HA) settings. Solvers ranked w.r.t. mean error.}
	\label{tab:prallc2}
% 	\customsize
	\small
% 	\footnotesize
	\centering
	\begin{tabular}{clcccc}
		\toprule
		Ranking & Solver & $\Bar{E}$[\%] & $\Bar{\sigma}$[\%] & $\Bar{\gamma}$[\%] & $\overline{CPU}$[$s$]  \\ 
		\midrule
21 & APSO & 1.173e+00 & 1.014e+00 & 45.9 & 1.0168 \\
2 & BARON (auto) & 2.054e-03 & 4.556e-17 & 92.0 & 1.7958 \\ 
5 & BARON (ipopt) & 2.055e-03 & 1.210e-09 & 92.0 & 1.9877 \\ 
3 & BARON (sd) & 2.054e-03 & 4.680e-17 & 92.0 & 1.8107 \\ 
4 & BARON (sqp) & 2.054e-03 & 4.120e-17 & 92.0 & 1.8730 \\ 
13 & FMINCON (active-set) & 1.770e-01 & 2.082e-01 & 72.3 & 0.0214 \\ 
10 & FMINCON (interior-point) & 1.985e-01 & 2.413e-01 & 75.9 & 0.0326 \\ 
11 & FMINCON (sqp) & 1.881e-01 & 2.388e-01 & 69.4 & 0.0082 \\ 
12 & FMINCON (sqp-legacy) & 1.857e-01 & 2.365e-01 & 69.7 & 0.0110 \\ 
17 & GCMMA & 5.112e-01 & 5.668e-01 & 45.2 & 1.0599 \\ 
8 & KNITRO (active-set) & 1.881e-01 & 2.691e-01 & 70.1 & 0.0698 \\ 
6 & KNITRO (interior-point) & 1.718e-01 & 1.962e-01 & 75.2 & 0.0357 \\ 
7 & KNITRO (sqp) & 1.785e-01 & 2.027e-01 & 73.0 & 0.1360 \\ 
15 & MIDACO & 2.735e-01 & 2.648e-01 & 85.9 & 1.0503 \\ 
20 & MMA & 9.786e-01 & 5.748e-01 & 42.1 & 0.7101 \\ 
18 & MQA & 5.358e-01 & 4.601e-01 & 20.8 & 3.2012 \\ 
1 & PENLAB & 1.502e-04 & 6.711e-39 & 31.0 & 0.0488 \\ 
19 & SGRA & 6.248e-01 & 7.673e-01 & 30.0 & 0.9632 \\ 
9 & SNOPT & 1.689e-01 & 2.010e-01 & 72.4 & 0.0069 \\ 
16 & SOLNP & 2.949e-01 & 3.106e-01 & 44.8 & 0.0112 \\ 
14 & SQA & 2.754e-01 & 2.838e-01 & 20.1 & 3.1838 \\ 
		\bottomrule    
	\end{tabular}
\end{table}
\FloatBarrier

\FloatBarrier
\begin{table}[ht]
	\caption{All constrained problems, quick solution (QS) settings. Solvers ranked w.r.t. mean CPU time.}
	\label{tab:prallc3}
% 	\customsize
	\small
% 	\footnotesize
	\centering
	\begin{tabular}{clcccc}
		\toprule
		Ranking & Solver & $\Bar{E}$[\%] & $\Bar{\sigma}$[\%] & $\Bar{\gamma}$[\%] & $\overline{CPU}$[$s$]  \\ 
		\midrule
10 & APSO & 1.531e+00 & 5.677e-01 & 35.2 & 0.0538 \\ 
17 & BARON (auto) & 7.925e-03 & 9.069e-05 & 92.0 & 0.7393 \\ 
16 & BARON (ipopt) & 1.152e-02 & 1.109e-03 & 92.0 & 0.7357 \\ 
18 & BARON (sd) & 2.766e-02 & 2.935e-04 & 92.0 & 0.7670 \\ 
19 & BARON (sqp) & 1.534e-02 & 7.491e-05 & 92.0 & 0.8652 \\ 
5 & FMINCON (active-set) & 2.850e-01 & 3.484e-01 & 68.9 & 0.0165 \\ 
6 & FMINCON (interior-point) & 2.166e-01 & 2.554e-01 & 72.3 & 0.0262 \\ 
2 & FMINCON (sqp) & 1.916e-01 & 2.448e-01 & 69.0 & 0.0071 \\ 
4 & FMINCON (sqp-legacy) & 1.902e-01 & 2.431e-01 & 69.1 & 0.0092 \\ 
15 & GCMMA & 6.967e-01 & 4.256e-01 & 45.5 & 0.5574 \\
8 & KNITRO (active-set) & 2.148e-01 & 2.767e-01 & 66.0 & 0.0295 \\
7 & KNITRO (interior-point) & 2.105e-01 & 2.826e-01 & 69.5 & 0.0268 \\ 
11 & KNITRO (sqp) & 2.207e-01 & 3.002e-01 & 70.3 & 0.0632 \\ 
12 & MIDACO & 8.355e-01 & 5.483e-01 & 58.6 & 0.0723 \\ 
14 & MMA & 1.161e+00 & 1.189e+00 & 41.0 & 0.1324 \\ 
20 & MQA & 5.844e-01 & 4.193e-01 & 20.8 & 3.1174 \\ 
9 & PENLAB & 1.896e-04 & 1.454e-37 & 31.0 & 0.0323 \\ 
13 & SGRA & 8.774e-01 & 1.198e+00 & 27.5 & 0.9369 \\ 
1 & SNOPT & 1.767e-01 & 2.045e-01 & 70.2 & 0.0027 \\ 
3 & SOLNP & 4.790e-01 & 6.452e-01 & 46.6 & 0.0087 \\ 
21 & SQA & 3.316e-01 & 3.131e-01 & 20.1 & 3.1361 \\
		\bottomrule    
	\end{tabular}
\end{table}
\FloatBarrier

% \begin{figure}[ht]
% 	\centering
% 	\includegraphics[width=\linewidth]{./figure/mean_error_c.png}
% 	\caption{Constrained problems, mean error.}
% 	\label{fig:meanEc}
% \end{figure}

% \begin{figure}[ht]
% 	\centering
% 	\includegraphics[width=\linewidth]{./figure/mean_var_c.png}
% 	\caption{Constrained problems, mean variance.}
% 	\label{fig:meanVc}
% \end{figure}

% \begin{figure}[ht]
% 	\centering
% 	\includegraphics[width=\linewidth]{./figure/mean_conv_c.png}
% 	\caption{Constrained problems, mean convergence rate.}
% 	\label{fig:meanCVc}
% \end{figure}

% \begin{figure}[ht]
% 	\centering
% 	\includegraphics[width=\linewidth]{./figure/mean_cpu_c.png}
% 	\caption{Constrained problems, mean CPU time.}
% 	\label{fig:meanCPUc}
% \end{figure}

% \clearpage

\section{Conclusions} 
\label{section:s6}
In this paper we provide an explicit comparison of a set of NLP solvers. The comparison includes popular solvers which are readily available in
MATLAB, a few gradient descent methods that have been extensively used in
literature, and a particle swarm optimization. Because of its widespread use
among research groups, both in academia and private sector, we have used MATLAB as common implementation platform. Constrained and unconstrained NLP problems have been selected amongst the standard benchmark problems with up to thirty variables and a up to nine scalar constraints. Results for the unconstrained problems show that BARON is the algorithm that deliver the best convergence rate and accuracy but it is the slowest.  PENLAB is the algorithm that has the best trade off between accuracy, convergence rate, and speed. For the constrained NLP problems, again, BARON is the solver which delivers excellent accuracy and convergence rate but is amongst the slower. FMINCON, KNITRO, SNOPT, and MIDACO are the one that are able to deliver a fair compromise of accuracy, convergence rate, and speed.

\section*{Data availability statement}
Data available on request from the authors.

\section*{Disclosure statement}
The authors declare that they have no known competing financial interests or personal relationships that could have appeared to influence the work reported in this paper. 

\section*{Funding}
This research received no external funding.

\bibliographystyle{tfcad}
\bibliography{mainNEW}% common bib file

\appendix

\section{Benchmark Test Functions for Unconstrained Global Optimization}
\label{sec:appendixA}

\subsection{Beale Function}
\begin{itemize}
    \item Dimension: $2$;
    \item Domain: $-4.5 \leq x_{i} \leq 4.5$;
    \item Function: 
    \begin{multline}
    f(\vec{x}) = \left( 1.5-x_{1}+x_{1}x_{2}\right)^2 + \left( 2.5-x_{1}+x_{1}x_{2}^2\right)^2 + \\ \left(2.625-x_{1}+x_{1}x_{2}^3\right)^2;
    \end{multline}
    \item Global minimum at $\vec{x^*} = (3, 0.5)$, $f(\vec{x^*}) = 0$.
\end{itemize}

\subsection{Bohachevsky 1 Function}
\begin{itemize}
    \item Dimension: $2$;
    \item Domain: $-100 \leq x_{i} \leq 100$;
    \item Function: 
    \begin{equation}
    f(\vec{x}) = x_{1}^2 + 2x_{2}^2 - 0.3\cos(3 \pi x_{1}) - 0.4\cos(4 \pi x_{2}) + 0.7;
    \end{equation}
    \item Global minimum at $\vec{x^*} = (0, 0)$, $f(\vec{x^*}) = 0$.
\end{itemize}

\subsection{Bohachevsky 2 Function}
\begin{itemize}
    \item Dimension: $2$;
    \item Domain: $-100 \leq x_{i} \leq 100$;
    \item Function: 
    \begin{equation}
    f(\vec{x}) = x_{1}^2 + 2x_{2}^2 - 0.3\cos(3 \pi x_{1})\cos(4 \pi x_{2}) + 0.3;
    \end{equation}
    \item Global minimum at $\vec{x^*} = (0, 0)$, $f(\vec{x^*}) = 0$.
\end{itemize}

\subsection{Bohachevsky 3 Function}
\begin{itemize}
    \item Dimension: $2$;
    \item Domain: $-100 \leq x_{i} \leq 100$;
    \item Function: 
    \begin{equation}
    f(\vec{x}) = x_{1}^2 + 2x_{2}^2 - 0.3\cos(3 \pi x_{1} + 4 \pi x_{2}) + 0.3;
    \end{equation}
    \item Global minimum at $\vec{x^*} = (0, 0)$, $f(\vec{x^*}) = 0$.
\end{itemize}

\subsection{Branin RCOS Function}
\begin{itemize}
    \item Dimension: $2$;
    \item Domain: $-5 \leq x_{1} \leq 10$, $0 \leq x_{2} \leq 15$;
    \item Function: 
    \begin{equation}
    f(\vec{x}) = {\left(-\frac{5.1\,{x_{1}}^2}{40\,\pi ^2}+\frac{5\,x_{1}}{\pi }+x_{2}-6\right)}^2+10\,\left(1-\frac{1}{8\,\pi }\right)\cos\left(x_{1}\right)+10;
    \end{equation}
    \item Global minimum at $\vec{x^*} = $ \{-pi, 12.275\}, \{pi, 2.275\}, \{9.42478,2.475\}, $f(\vec{x^*}) = $ 0.397887.
\end{itemize}

\subsection{Colville Function}
\begin{itemize}
    \item Dimension: $4$;
    \item Domain: $-10 \leq x_{i} \leq 10$;
    \item Function: 
    \begin{multline}
    f(\vec{x}) = 100(x_{1}^2-x_{2})^2+(x_{1}-1)^2+(x_{3}-1)^2+\\90(x_{3}^2-x_{4})^2+10.1((x_{2}-1)^2+(x_{4}-1)^2)+19.8(x_{2}-1)(x_{4}-1);
    \end{multline}
    \item Global minimum at $\vec{x^*} = (1, 1, 1, 1)$, $f(\vec{x^*}) = 0$.
\end{itemize}

\subsection{Dixon \& Price Function}
\begin{itemize}
    \item Dimension: $25$;
    \item Domain: $-10 \leq x_{i} \leq 10$;
    \item Function: 
    \begin{equation}
    f(\vec{x}) = (x_{1}-2)^2 + \sum_{i=1}^{n}i(2x_{i}^2-x_{i-1})^2;
    \end{equation}
    \item Global minimum at $\vec{x^*} = 2^{-\frac{2^i-2}{2^i}}$ with $i=1,...,n$, $f(\vec{x^*}) = 0$.
\end{itemize}

\subsection{Hump Function}
\begin{itemize}
    \item Dimension: $2$;
    \item Domain: $-5 \leq x_{i} \leq 5$;
    \item Function: 
    \begin{equation}
    f(\vec{x}) = x_{1}^2(x_{1}^{\frac{4}{3}} - 2.1x_{1}^2 + 4) + x_{1}x_{2} + x_{2}^2(4x_{2}^2 - 4) ;
    \end{equation}
    \item Global minimum at $\vec{x^*} = \{0.0898, -0.7126\}, \{-0.0898, 0.7126\}$, $f(\vec{x^*}) = 0$.
\end{itemize}

\subsection{Matyas Function}
\begin{itemize}
    \item Dimension: $2$;
    \item Domain: $-10 \leq x_{i} \leq 10$;
    \item Function: 
    \begin{equation}
    f(\vec{x}) = 0.26(x_{1}^2+x_{2}^2)-0.48x_{1}x_{2};
    \end{equation}
    \item Global minimum at $\vec{x^*} = (0, 0)$, $f(\vec{x^*}) = 0$.
\end{itemize}

\subsection{Perm(n,\textbeta) Function}
\begin{itemize}
    \item Dimension: $10$;
    \item Domain: $-10 \leq x_{i} \leq 10$;
    \item Function: 
    \begin{equation}
    f(\vec{x}) = \sum_{k=1}^n \left(\sum_{i=1}^n (i^k+\beta)\left( \frac{x_i}{i} ^k-1\right)\right)^2, \quad with \quad \beta = 0.5;
    \end{equation}
    \item Global minimum at $\vec{x^*} = i$ with $i=1,...,n$, $f(\vec{x^*}) = 0$.
\end{itemize}

\subsection{Powell singular Function}
\begin{itemize}
    \item Dimension: $16$;
    \item Domain: $-4 \leq x_{i} \leq 5$;
    \item Function: 
    \begin{multline}
    f(\vec{x}) = \sum_{i=1}^{n/4} (x_{4i-3} + 10x_{4i-2})^2 + 5(x_{4i-1}-x_{4i})^2 + \\ (x_{4i-2}-x_{4i-1})^2 + 10(x_{4i-3}-x_{4i})^2;
    \end{multline}
    \item Global minimum at $\vec{x^*} = (0, ..., 0)$, $f(\vec{x^*}) = 0$.
\end{itemize}

\subsection{Power Sum Function}
\begin{itemize}
    \item Dimension: $4$;
    \item Domain: $0 \leq x_{i} \leq 256$;
    \item Function: 
    \begin{equation}
    f(\vec{x}) = \sum_{k=1}^n \left(\left(\sum_{i=1}^n x_i^k\right) - b_k\right)^2, \quad with \quad b = (8,18,44,114);
    \end{equation}
    \item Global minimum at $\vec{x^*} = (1, 2, 3, 3)$, $f(\vec{x^*}) = 0$.
\end{itemize}

\subsection{Sphere Function}
\begin{itemize}
    \item Dimension: $30$;
    \item Domain: $-5.12 \leq x_{i} \leq 5.12$;
    \item Function: 
    \begin{equation}
    f(\vec{x}) = \sum_{i=1}^n x_{i}^2;
    \end{equation}
    \item Global minimum at $\vec{x^*} = (0, ..., 0)$, $f(\vec{x^*}) = 0$.
\end{itemize}

\subsection{Sum Squares Function}
\begin{itemize}
    \item Dimension: $30$;
    \item Domain: $-10 \leq x_{i} \leq 10$;
    \item Function: 
    \begin{equation}
    f(\vec{x}) = \sum_{i=1}^n ix_{i}^2;
    \end{equation}
    \item Global minimum at $\vec{x^*} = (0, ..., 0)$, $f(\vec{x^*}) = 0$.
\end{itemize}

\subsection{Trid Function}
\begin{itemize}
    \item Dimension: $10$;
    \item Domain: $-100 \leq x_{i} \leq 100$;
    \item Function: 
    \begin{equation}
    f(\vec{x}) = \sum_{i=1}^n (x_i-1)^2 - \sum_{i=2}^n x_ix_{i-1};
    \end{equation}
    \item Global minimum at $\vec{x^*} = i*\left(11-i\right)$ with $i = 1, ..., n$, $f(\vec{x^*}) = -210$.
\end{itemize}

\subsection{Zakharov Function}
\begin{itemize}
    \item Dimension: $20$;
    \item Domain: $-5 \leq x_{i} \leq 10$;
    \item Function: 
    \begin{equation}
    f(\vec{x}) = \sum_{i=1}^n x_{i}^2 + \left(\frac{1}{2}\sum_{i=1}^n ix_{i}^2\right)^2 + \left(\frac{1}{2}\sum_{i=1}^n ix_{i}^2\right)^4;
    \end{equation}
    \item Global minimum at $\vec{x^*} = (0, ..., 0)$, $f(\vec{x^*}) = 0$.
\end{itemize}

\subsection{Branin RCOS 2 Function}
\begin{itemize}
    \item Dimension: $2$;
    \item Domain: $-5 \leq x_{i} \leq 15$;
    \item Function: 
    \begin{multline}
    f(\vec{x}) = {\left(-\frac{5.1\,{x_{1}}^2}{40\,\pi ^2}+\frac{5\,x_{1}}{\pi }+x_{2}-6\right)}^2+\\10\,\left(1-\frac{1}{8\,\pi }\right)\cos\left(x_{1}\right)\cos\left(x_{2}\right)\ln(x_{1}^2+x_{2}^2+1)+10;
    \end{multline}
    \item Global minimum at $\vec{x^*} = (-3.2, 12.53)$, $f(\vec{x^*}) = 5.559037$.
\end{itemize}

\subsection{Ackley 1 Function}
\begin{itemize}
    \item Dimension: $10$;
    \item Domain: $-15 \leq x_{i} \leq 30$;
    \item Function: 
    \begin{equation}
    f(\vec{x}) = -20 e^{-0.2 \sqrt{\frac{1}{n}\sum_{i=0}^{n}x_i^{2}}} - e^{\frac{1}{n}\sum_{i=0}^{n} \cos(2\pi x_i)} + 20+e;
    \end{equation}
    \item Global minimum at $\vec{x^*} = (0, ..., 0)$, $f(\vec{x^*}) = 0$.
\end{itemize}

\subsection{Ackley 2 Function}
\begin{itemize}
    \item Dimension: $2$;
    \item Domain: $-32 \leq x_{i} \leq 32$;
    \item Function: 
    \begin{equation}
    f(\vec{x}) = -200 e^{-0.02 \sqrt{x_1^{2}+x_2^{2}}};
    \end{equation}
    \item Global minimum at $\vec{x^*} = (0, 0)$, $f(\vec{x^*}) = -200$.
\end{itemize}

\subsection{Camel 3 Function}
\begin{itemize}
    \item Dimension: $2$;
    \item Domain: $-5 \leq x_{i} \leq 5$;
    \item Function: 
    \begin{equation}
    f(\vec{x}) = 2x_1^{2} -1.05x_1^{4} + \frac{1}{6}x_1^{6} + x_1x_2 + x_2^{2};
    \end{equation}
    \item Global minimum at $\vec{x^*} = (0, 0)$, $f(\vec{x^*}) = 0$.
\end{itemize}

\subsection{Booth Function}
\begin{itemize}
    \item Dimension: $2$;
    \item Domain: $-10 \leq x_{i} \leq 10$;
    \item Function: 
    \begin{equation}
    f(\vec{x}) = \left( x_1 + 2x_2 -7 \right)^2 + \left( 2x_1 + x_2 -5 \right)^2;
    \end{equation}
    \item Global minimum at $\vec{x^*} = (1, 3)$, $f(\vec{x^*}) = 0$.
\end{itemize}

\subsection{Brown Function}
\begin{itemize}
    \item Dimension: $14$;
    \item Domain: $-1 \leq x_{i} \leq 4$;
    \item Function: 
    \begin{equation}
    f(\vec{x}) = \sum_{i=0}^{n-1} \left( x_i^2\right)^{\left( x_{i+1}^2+1\right)} + \left( x_{i+1}^2\right)^{\left( x_{i}^2+1\right)} ;
    \end{equation}
    \item Global minimum at $\vec{x^*} = (0, ..., 0)$, $f(\vec{x^*}) = 0$.
\end{itemize}

\subsection{Cube Function}
\begin{itemize}
    \item Dimension: $2$;
    \item Domain: $-10 \leq x_{i} \leq 10$;
    \item Function: 
    \begin{equation}
    f(\vec{x}) = 100\left( x_2 - x_1^3 \right)^2 + \left( 1 - x_1\right)^2;
    \end{equation}
    \item Global minimum at $\vec{x^*} = (-1, 1)$, $f(\vec{x^*}) = 0$.
\end{itemize}

\subsection{Exponential Function}
\begin{itemize}
    \item Dimension: $18$;
    \item Domain: $-1 \leq x_{i} \leq 1$;
    \item Function: 
    \begin{equation}
    f(\vec{x}) = -e^{\left( -0.5 \sum_{i=1}^{n} x_i^2 \right)};
    \end{equation}
    \item Global minimum at $\vec{x^*} = (0, ..., 0)$, $f(\vec{x^*}) = 1$.
\end{itemize}

\subsection{Freudenstein Roth Function}
\begin{itemize}
    \item Dimension: $2$;
    \item Domain: $-10 \leq x_{i} \leq 10$;
    \item Function: 
    \begin{multline}
    f(\vec{x}) = \left( x_1 - 13 + x_2\left(\left( 5-x_2\right)x_2 -2\right) \right)^2 + \\ \left( x_1 - 29 + x_2\left(\left( x_2+1\right)x_2 -14\right) \right)^2;
    \end{multline}
    \item Global minimum at $\vec{x^*} = (5, 4)$, $f(\vec{x^*}) = 0$.
\end{itemize}

\subsection{Miele Cantrell Function}
\begin{itemize}
    \item Dimension: $4$;
    \item Domain: $-1 \leq x_{i} \leq 1$;
    \item Function: 
    \begin{equation}
    f(\vec{x}) = \left(e^{-x_1} -x_2\right)^4 +100\left(x_2 -x_3\right)^6 + \left(\tan(x_3 -x_4)\right)^4 +x_1^{8};
    \end{equation}
    \item Global minimum at $\vec{x^*} = (0, 1, 1, 1)$, $f(\vec{x^*}) = 0$.
\end{itemize}

\subsection{Quadratic Function}
\begin{itemize}
    \item Dimension: $2$;
    \item Domain: $-10 \leq x_{i} \leq 10$;
    \item Function: 
    \begin{multline}
    f(\vec{x}) = - 3803.84 - 138.08x_1 - 232.92x_2 + \\ 128.08x_1^2 + 203.64x_2^2 + 182.25x_1x_2;
    \end{multline}
    \item Global minimum at $\vec{x^*} = (0.19388, 0.48513)$, $f(\vec{x^*}) = -3873.7243$.
\end{itemize}

\subsection{Rotated Ellipse Function}
\begin{itemize}
    \item Dimension: $2$;
    \item Domain: $-500 \leq x_{i} \leq 500$;
    \item Function: 
    \begin{equation}
    f(\vec{x}) = 7x_1^2 - 6\sqrt{3}x_1x_2 + 13x_2^2;
    \end{equation}
    \item Global minimum at $\vec{x^*} = (0, 0)$, $f(\vec{x^*}) = 0$.
\end{itemize}

\subsection{Rump Function}
\begin{itemize}
    \item Dimension: $2$;
    \item Domain: $-500 \leq x_{i} \leq 500$;
    \item Function: 
    \begin{equation}
    f(\vec{x}) = (333.75-x_1^2)x_2^6 + x_1^2(11x_1^2x_2^2 -121x_2^4 -2) + 5.5x_2^8 + \frac{x_1}{2+x_2};
    \end{equation}
    \item Global minimum at $\vec{x^*} = (0, 0)$, $f(\vec{x^*}) = 0$.
\end{itemize}

\subsection{Wayburn Seader 3 Function}
\begin{itemize}
    \item Dimension: $2$;
    \item Domain: $-500 \leq x_{i} \leq 500$;
    \item Function: 
    \begin{equation}
    f(\vec{x}) = 2\frac{x_1^3}{3} - 8x_1^2 + 33x_1 - x_1x_2 + 5 + \left[ (x_1-4)^2 + (x_2-5)^2 -4 \right]^2;
    \end{equation}
    \item Global minimum at $\vec{x^*} = f(5.611, 6.187)$, $f(\vec{x^*}) = 21.35$.
\end{itemize}

\section{Benchmark Test Functions for Constrained Global Optimization}
\label{sec:appendixB}

% \subsection{G1 Problem}
\subsection{}
\begin{itemize}
    \item Dimension: $13$;
    \item Domain and search space: $0 \leq x_{i} \leq u_{i}$, with \\ $u=(1,1,1,...,1,100,100,100,1)$;
    \item Function: 
    \begin{equation}
    f(\vec{x}) = 5\sum_{i=1}^{4}x_{i} - 5\sum_{i=1}^{4}x_{i}^2 - \sum_{i=5}^{n}x_{i};
    \end{equation}
    \item Constraints: 
    \begin{equation}\begin{cases}
    c_{1}(\vec{x}) = 2x_{1}+2x_{2}+x_{10}+x_{11}-10 \leq 0; \\
    c_{2}(\vec{x}) = 2x_{1}+2x_{3}+x_{10}+x_{12}-10 \leq 0; \\
    c_{3}(\vec{x}) = 2x_{2}+2x_{3}+x_{11}+x_{12}-10 \leq 0; \\
    c_{4}(\vec{x}) = -8x_{1}+x_{10} \leq 0; \\
    c_{5}(\vec{x}) = -8x_{2}+x_{11} \leq 0; \\
    c_{6}(\vec{x}) = -8x_{3}+x_{12} \leq 0; \\
    c_{7}(\vec{x}) = -2x_{4}-x_{5}+x_{10} \leq 0; \\
    c_{8}(\vec{x}) = -2x_{6}-x_{7}+x_{11} \leq 0; \\
    c_{9}(\vec{x}) = -2x_{8}-x_{9}+x_{12} \leq 0;
    \end{cases}\end{equation}
    \item Global minimum at $\vec{x^*} = (1,1,1,1,1,1,1,1,1,3,3,3,1)$, $f(\vec{x^*}) = -15$.
\end{itemize}

% \subsection{G2 Problem}
\subsection{}
\begin{itemize}
    \item Dimension: $20$;
    \item Domain and search space: $0 \leq x_{i} \leq 10$;
    \item Function: 
    \begin{equation}
    f(\vec{x}) = -\abs{\frac{\sum_{i=1}^{n}\cos(x_{i})^4-2\prod_{i=1}^{n}\cos(x_{i})^2}{\sqrt{\sum_{i=1}^{n}ix_{i}^2}}};
    \end{equation}
    \item Constraints: 
    \begin{equation}\begin{cases}
    c_{1}(\vec{x}) = -\prod_{i=1}^{n}x_{i} + 0.75 \leq 0; \\
    c_{2}(\vec{x}) = \sum_{i=1}^{n}x_{i} - 7.5n \leq 0;
    \end{cases}\end{equation}
    \item Global minimum (best known) at $f(\vec{x^*}) = -0.803619$.
\end{itemize}

% \subsection{G3 Problem}
\subsection{}
\begin{itemize}
    \item Dimension: $10$;
    \item Domain and search space: $0 \leq x_{i} \leq 1$;
    \item Function: 
    \begin{equation}
    f(\vec{x}) = -\left(\sqrt{n}\right)^n \prod_{i=1}^{n}x_{i};
    \end{equation}
    \item Constraints: 
    \begin{equation}
    c_{1}(\vec{x}) = \sum_{i=1}^{n}x_{i}^2-1 = 0;
    \end{equation}
    \item Global minimum at $\vec{x^*} = \left(\left(\frac{1}{10}\right)^{0.5},...,\left(\frac{1}{10}\right)^{0.5}\right)$, $f(\vec{x^*}) = -1$.
\end{itemize}

% \subsection{G4 Problem}
\subsection{}
\begin{itemize}
    \item Dimension: $5$;
    \item Domain and search space: $l_{i} \leq x_{i} \leq u_{i}$, with $l=(78,33,27,27,27)$ \\ and $u=(102,45,45,45,45)$;
    \item Function: 
    \begin{equation}
    f(\vec{x}) = 5.3578547x_{3}^2+0.8356891x_{1}x_{5}+37.293239x_{1}-40792.141;
    \end{equation}
    \item Constraints: 
    \begin{equation}\begin{cases}
    c_{1}(\vec{x}) = -u \leq 0; \\
    c_{2}(\vec{x}) = u-92 \leq 0; \\
    c_{3}(\vec{x}) = -v+90 \leq 0; \\
    c_{4}(\vec{x}) = v-110 \leq 0; \\
    c_{5}(\vec{x}) = -w+20 \leq 0; \\
    c_{6}(\vec{x}) = w-25 \leq 0; \\
    u = 85.334407+0.0056858x_{2}x_{5}+0.0006262x_{1}x_{4}-0.0022053x_{3}x_{5}; \\
    v = 80.51249+0.0071317x_{2}x_{5}+0.0029955x_{1}x_{2}+0.0021813x_{3}^2; \\
    w = 9.300961+0.0047026x_{3}x_{5}+0.0012547x_{1}x_{3}+0.0019085x_{3}x_{4};
    \end{cases}\end{equation}
    \item Global minimum at $\vec{x^*} = (78,33,29.995,45,36.7758)$, $f(\vec{x^*}) =  -30665.539$.
\end{itemize}

% \subsection{G5 Problem}
\subsection{}
\begin{itemize}
    \item Dimension: $4$;
    \item Domain and search space: $l_{i} \leq x_{i} \leq u_{i}$, with $l=(0,0,-0.55,-0.55)$ \\ and $u=(1200,1200,0.55,0.55)$;
    \item Function: 
    \begin{equation}
    f(\vec{x}) = 3x_{1}+10^{-6}x_{1}^3+2x_{2}+\frac{2}{3}10^{-6}x_{2}^3;
    \end{equation}
    \item Constraints: 
    \begin{equation}\begin{cases}
    c_{1}(\vec{x}) = x_{3} - x_{4} - 0.55 \leq 0; \\
    c_{2}(\vec{x}) = x_{4} - x_{3} - 0.55 \leq 0; \\
    c_{3}(\vec{x}) = 1000(\sin(-x_{3}-0.25)+\sin(-x_{4}-0.25))+894.8-x_{1} = 0; \\
    c_{4}(\vec{x}) = 1000(\sin(x_{3}-0.25)+\sin(x_{3}-x_{4}-0.25))+894.8-x_{2} = 0; \\
    c_{5}(\vec{x}) = 1000(\sin(x_{4}-0.25)+\sin(x_{4}-x_{3}-0.25))+1294.8 = 0;
    % c_{3}(\vec{x}) = \abs{1000(\sin(-x_{3}-0.25)+\sin(-x_{4}-0.25))+894.8-x_{1}}-10^{-4} = 0; \\
    % c_{4}(\vec{x}) = \abs{1000(\sin(x_{3}-0.25)+\sin(x_{3}-x_{4}-0.25))+894.8-x_{2}}-10^{-4} = 0; \\
    % c_{5}(\vec{x}) = \abs{1000(\sin(x_{4}-0.25)+\sin(x_{4}-x_{3}-0.25))+1294.8}-10^{-4} = 0;
    \end{cases}\end{equation}
    \item Global minimum at $\vec{x^*} = (679.9453,1026,0.118876,-0.3962336)$, $f(\vec{x^*}) = 5126.4981$.
\end{itemize}

% \subsection{G6 Problem}
\subsection{}
\begin{itemize}
    \item Dimension: $2$;
    \item Domain and search space: $l_{i} \leq x_{i} \leq 100$, with $l=(13,0)$;
    \item Function: 
    \begin{equation}
    f(\vec{x}) = (x_{1}-10)^3 + (x_{2}-20)^3;
    \end{equation}
    \item Constraints: 
    \begin{equation}\begin{cases}
    c_{1}(\vec{x}) = -(x_{1}-5)^2-(x_{2}-5)^2+100 \leq 0; \\
    c_{2}(\vec{x}) = (x_{1}-6)^2+(x_{2}-5)^2-82.81 \leq 0;
    \end{cases}\end{equation}
    \item Global minimum at $\vec{x^*} = (14.095,0.84296)$, $f(\vec{x^*}) = -6961.81388$.
\end{itemize}

% \subsection{G7 Problem}
\subsection{}
\begin{itemize}
    \item Dimension: $10$;
    \item Domain and search space: $-10 \leq x_{i} \leq 10$;
    \item Function: 
    \begin{multline}
    f(\vec{x}) = 2(x_{6} - 1)^2 - 16x_{2} - 14x_{1} + (x_{5} - 3)^2 + 4(x_{4} - 5)^2 + (x_{3} - 10)^2 + \\ (x_{10} - 7)^2 + 7(x_{8} - 11)^2 + 2(x_{9} - 10)^2 + x_{1}x_{2} + x_{1}^2 + x_{2}^2 + 5x_{7}^2 + 45;
    \end{multline}
    \item Constraints: 
    \begin{equation}\begin{cases}
    c_{1}(\vec{x}) = 4x_{1}+5x_{2}-3x_{7}+9x_{8}-105 \leq 0; \\
    c_{2}(\vec{x}) = 10x_{1}-8x_{2}-17x_{7}+2x_{8} \leq 0; \\  
    c_{3}(\vec{x}) = -8x_{1}+2x_{2}+5x_{9}-2x_{10}-12 \leq 0; \\ 
    c_{4}(\vec{x}) = 3(x_{1}-2)^2+4(x_{2}-3)^2+2x_{3}^2-7x_{4}-120 \leq 0; \\      
    c_{5}(\vec{x}) = 5x_{1}^2+8x_{2}+(x_{3}-6)^2-2x_{4}-40 \leq 0; \\
    c_{6}(\vec{x}) = 0.5(x_{1}-8)^2+2(x_{2}-4)^2+3x_{5}^2-x_{6}-30 \leq 0; \\ 
    c_{7}(\vec{x}) = x_{1}^2+2(x_{2}-2)^2-2x_{1}x_{2}+14x_{5}-6x_{6} \leq 0; \\       
    c_{8}(\vec{x}) = -3x_{1}+6x_{2}+12(x_{9}-8)^2-7x_{10} \leq 0;
    \end{cases}\end{equation}
    \item Global minimum at $\vec{x^*} =$ (2.171996, 2.363683, 8.773926, 5.095984, 0.9906548, 1.430574, 1.321644, 9.828726, 8.280092, 8.375927), $f(\vec{x^*}) =$ 24.3062091.
\end{itemize}

% \subsection{G8 Problem}
\subsection{}
\begin{itemize}
    \item Dimension: $2$;
    \item Domain and search space: $0 \leq x_{i} \leq 10$;
    \item Function: 
    \begin{equation}
    f(\vec{x}) = -\frac{\sin(2\pi x_{1})^3\sin(2\pi x_{2})}{x_{1}^3(x_{1}+x_{2})};
    \end{equation}
    \item Constraints: 
    \begin{equation}\begin{cases}
    c_{1}(\vec{x}) = x_{1}^2 - x_{2} + 1 \leq 0; \\
    c_{2}(\vec{x}) = 1- x_{1} + (x_{2}-4)^2 \leq 0;
    \end{cases}\end{equation}
    \item Global minimum at $\vec{x^*} = (1.2279713, 4.2453733)$, $f(\vec{x^*}) = -0.095825$.
\end{itemize}

% \subsection{G9 Problem}
\subsection{}
\begin{itemize}
    \item Dimension: $7$;
    \item Domain and search space: $-10 \leq x_{i} \leq 10$;
    \item Function: 
    \begin{multline}
    f(\vec{x}) = (x_{1} - 10)^2 - 8x_{7} - 10x_{6} + 5(x_{2} - 12)^2 + \\ 3(x_{4} - 11)^2 - 4x_{6}x_{7} + x_{3}^4 + 7x_{6}^2 + 10x_{5}^6 + x_{7}^4;
    \end{multline}
    \item Constraints: 
    \begin{equation}\begin{cases}
    c_{1}(\vec{x}) = 2x_{1}^2+3x_{2}^4+x_{3}+4x_{4}^2+5x_{5}-127 \leq 0; \\
    c_{2}(\vec{x}) = 7x_{1}+3x_{2}+10x_{3}^2+x_{4}-x_{5}-282 \leq 0; \\
    c_{3}(\vec{x}) = 23x_{1}+x_{2}^2+6x_{6}^2-8x_{7}-196 \leq 0; \\
    c_{4}(\vec{x}) = 4x_{1}^2+x_{2}^2-3x_{1}x_{2}+2x_{3}^2+5x_{6}-11x_{7} \leq 0;
    \end{cases}\end{equation}
    \item Global minimum at $\vec{x^*} =$ (2.330499, 1.951372, -0.4775414,  4.365726, \\ -0.6244870, 1.038131, 1.594227), $f(\vec{x^*}) =$ 680.6300573.
\end{itemize}

% \subsection{G10 Problem}
\subsection{}
\begin{itemize}
    \item Dimension: $8$;
    \item Domain and search space: $l_{i} \leq x_{i} \leq u_{i}$, with $l=10(10,100,100,1,1,1,1,1)$ and $u=1000(10,10,10,1,1,1,1,1)$;
    \item Function: 
    \begin{equation}
    f(\vec{x}) = x_{1} + x_{2} + x_{3};
    \end{equation}
    \item Constraints: 
    \begin{equation}\begin{cases}
    c_{1}(\vec{x}) = -1+0.0025(x_{4}+x_{6}) \leq 0; \\
    c_{2}(\vec{x}) = -1+0.0025(-x_{4}+x_{5}+x_{7}) \leq 0; \\
    c_{3}(\vec{x}) = -1+0.01(-x_{5}+x_{8}) \leq 0; \\
    c_{4}(\vec{x}) = 100x_{1}-x_{1}x_{6}+833.33252x_{4}-83333.333 \leq 0; \\
    c_{5}(\vec{x}) = x_{2}x_{4}-x_{2}x_{7}-1250x_{4}+1250x_{5} \leq 0; \\
    c_{6}(\vec{x}) = x_{3}x_{5}-x_{3}x_{8}-2500x_{5}+1250000 \leq 0;
    \end{cases}\end{equation}
    \item Global minimum at $\vec{x^*} =$ (579.3167, 1359.943, 5110.071, 182.0174, 295.5985, 217.9799, 286.4162,395.5979), $f(\vec{x^*}) = 7049.3307$.
\end{itemize}

% \subsection{G11 Problem}
\subsection{}
\begin{itemize}
    \item Dimension: $2$;
    \item Domain and search space: $-1 \leq x_{i} \leq 1$;
    \item Function: 
    \begin{equation}
    f(\vec{x}) = x_{1}^2 + (x_{2}-1)^2;
    \end{equation}
    \item Constraints: 
    \begin{equation}
    c_{1}(\vec{x}) = x_{2} - x_{1}^2 = 0;
    \end{equation}
    \item Global minimum at $\vec{x^*} = \pm(\frac{1}{2}^{0.5}, \frac{1}{2})$, $f(\vec{x^*}) = 0.75$.
\end{itemize}

% \subsection{G13 Problem}
\subsection{}
\begin{itemize}
    \item Dimension: $5$;
    \item Domain and search space: $l_{i} \leq x_{i} \leq u_{i}$, with $l=-u$, $u=$(2.3,2.3,3.2,3.2,3.2);
    \item Function: 
    \begin{equation}
    f(\vec{x}) = e^{x_{1}x_{2}x_{3}x_{4}x_{5}};
    \end{equation}
    \item Constraints: 
    \begin{equation}\begin{cases}
    c_{1}(\vec{x}) = x_{1}^2 + x_{2}^2 + x_{3}^2 + x_{4}^2 + x_{5}^2 -10 = 0; \\
    c_{2}(\vec{x}) = x_{2}x_{3} -5x_{4}x_{5} = 0; \\
    c_{3}(\vec{x}) = x_{1}^3 + x_{2}^3 + 1 = 0;
    \end{cases}\end{equation}
    \item Global minimum at $\vec{x^*} =$ (-1.717143, 1.595709, 1.827247, -0.7636413, \\ -0.763645), $f(\vec{x^*}) = 0.0539498$.
\end{itemize}

% \subsection{G1-Ciarcia Problem}
\subsection{}
\begin{itemize}
    \item Dimension: $2$;
    \item Domain and search space: $-10 \leq x_{i} \leq 10$;
    \item Function: 
    \begin{equation}
    f(\vec{x}) = x_{1}^3 + x_{2}^3;
    \end{equation}
    \item Constraints: 
    \begin{equation}
    c_{1}(\vec{x}) = x_{1} + x_{2} - 8 = 0;
    \end{equation}
    \item Global minimum at $\vec{x^*} =$ (4, 4), $f(\vec{x^*}) = 128$.
\end{itemize}

% \subsection{G2-Ciarcia Problem}
\subsection{}
\begin{itemize}
    \item Dimension: $3$;
    \item Domain and search space: $-10 \leq x_{i} \leq 10$;
    \item Function: 
    \begin{equation}
    f(\vec{x}) = (x_{1} - 1)^2 + (x_{2} - 2)^2 + x_{3}^2 +2;
    \end{equation}
    \item Constraints: 
    \begin{equation}
    c_{1}(\vec{x}) = x_{1}^2 + x_{2} - 3 = 0;
    \end{equation}
    \item Global minimum at $\vec{x^*} =$ (1, 2, 0), $f(\vec{x^*}) = 2$.
\end{itemize}

% \subsection{G3-Ciarcia Problem}
\subsection{}
\begin{itemize}
    \item Dimension: $3$;
    \item Domain and search space: $-10 \leq x_{i} \leq 10$;
    \item Function: 
    \begin{equation}
    f(\vec{x}) = 2(x_{1}x_{2} + x_{2}x_{3} + x_{1}x_{3});
    \end{equation}
    \item Constraints: 
    \begin{equation}\begin{cases}
    c_{1}(\vec{x}) = x_{1}x_{2}x_{3} - 72 = 0; \\
    c_{2}(\vec{x}) = x_{1} - 2x_{2} = 0;
    \end{cases}\end{equation}
    \item Global minimum at $\vec{x^*} =$ (6, 3, 4), $f(\vec{x^*}) = 108$.
\end{itemize}

% \subsection{Schittkowski-7 Problem}
\subsection{}
\begin{itemize}
    \item Dimension: $2$;
    \item Domain and search space: $-10 \leq x_{i} \leq 10$;
    \item Function: 
    \begin{equation}
    f(\vec{x}) = \ln \left( 1 + x_{1}^2 \right) - x_{2};
    \end{equation}
    \item Constraints: 
    \begin{equation}
    c_{1}(\vec{x}) = (1+x_{1}^2)^2 + x_{2}^2 -4  = 0
    \end{equation}
    \item Global minimum at $\vec{x^*} = (0, \sqrt{3})$, $f(\vec{x^*}) = -\sqrt{3}$.
\end{itemize}

% \subsection{Schittkowski-27 Problem}
\subsection{}
\begin{itemize}
    \item Dimension: $3$;
    \item Domain and search space: $-10 \leq x_{i} \leq 10$;
    \item Function: 
    \begin{equation}
    f(\vec{x}) = 0.01(x_{1}-1)^2 + (x_{2}-x_{1}^2)^2;
    \end{equation}
    \item Constraints: 
    \begin{equation}
    c_{1}(\vec{x}) = x_{1} + x_{3}^2 + 1  = 0
    \end{equation}
    \item Global minimum at $\vec{x^*} =$ (-1, 1, 0), $f(\vec{x^*}) =$ 0.04.
\end{itemize}

% \subsection{2.2-TP1 Problem}
\subsection{}
\begin{itemize}
    \item Dimension: $5$;
    \item Domain and search space: $0 \leq x_{i} \leq 1$;
    \item Function: 
    \begin{equation}
    f(\vec{x}) = \vec{c}^T\vec{x} - 0.5\vec{x}^T\vec{Q}\vec{x};
    \end{equation}
    \begin{equation*}
    \text{with} \quad \vec{c} = (42, 44, 45, 47, 47.5)^T, \quad \vec{Q}=100\vec{I}
    \end{equation*}
    \item Constraints: 
    \begin{equation}
    c_{1}(\vec{x}) = 20x_{1} + 12x_{2} + 11x_{3} + 7x_{4} + 4x_{5} - 40 \leq 0
    \end{equation}
    \item Global minimum at $\vec{x^*} =$ (1, 1, 0, 1, 0), $f(\vec{x^*}) = -17$.
\end{itemize}

% \subsection{2.3-TP2 Problem}
\subsection{}
\begin{itemize}
    \item Dimension: $6$;
    \item Domain and search space: $0 \leq x_{1,\dots,5} \leq 1$, $0 \leq y$;
    \item Function: 
    \begin{equation}
    f(\vec{x},y) = \vec{c}^T\vec{x} - 0.5\vec{x}^T\vec{Q}\vec{x} -10y;
    \end{equation}
    \begin{equation*}
    \text{with} \quad \vec{c} = (-10.5, -7.5, -3.5, -2.5, -1.5)^T, \quad \vec{Q}=\vec{I}
    \end{equation*}
    \item Constraints: 
    \begin{equation}\begin{cases}
    c_{1}(\vec{x}) = 6x_{1} + 3x_{2} + 3x_{3} + 2x_{4} + 1x_{5} - 6.5 \leq 0 \\
    c_{2}(\vec{x}) = 10x_{1} + 10x_{3} + y - 20 \leq 0
    \end{cases}\end{equation}
    \item Global minimum at $\vec{x^*} =$ (0, 1, 0, 1, 1), $y^* =20$, $f(\vec{x^*}) = -213$.
\end{itemize}

% \subsection{2.4-TP3 Problem}
\subsection{}
\begin{itemize}
    \item Dimension: $12$;
    \item Domain and search space: $0 \leq x_{i} \leq 1$, $0 \leq y_{1,\dots,5} \leq 1$, $0 \leq y_{6,\dots,8} \leq 3$;
    \item Function: 
    \begin{equation}
    f(\vec{x},y) = \vec{c}^T\vec{x} - 0.5\vec{x}^T\vec{Q}\vec{x} +\vec{d}^T\vec{y};
    \end{equation}
    \begin{equation*}
    \text{with} \quad \vec{c} = (5,5,5,5)^T, \quad \vec{Q}=100\vec{I}, 
    \end{equation*}
    \begin{equation*}
    \vec{d}=(-1,-1,-1,-1,-1,-1,-1,-1)^T
    \end{equation*}
    \item Constraints: 
    \begin{equation}\begin{cases}
    c_{1}(\vec{x}) = 2x_{1} + 2x_{2} + y_{6} + y_{7} - 10 \leq 0 \\
    c_{2}(\vec{x}) = 2x_{1} + 2x_{3} + y_{6} + y_{8} - 10 \leq 0 \\
    c_{3}(\vec{x}) = 2x_{2} + 2x_{3} + y_{7} + y_{8} - 10 \leq 0 \\
    c_{4}(\vec{x}) = -8x_{1} + y_{6} \leq 0 \\
    c_{5}(\vec{x}) = -8x_{2} + y_{7} \leq 0 \\
    c_{6}(\vec{x}) = -8x_{3} + y_{8} \leq 0 \\
    c_{7}(\vec{x}) = -2x_{4} - y_{1} + y_{6} \leq 0 \\
    c_{8}(\vec{x}) = -2y_{2} - y_{3} + y_{7} \leq 0 \\
    c_{9}(\vec{x}) = -2y_{4} - y_{5} + y_{8} \leq 0 \\
    \end{cases}\end{equation}
    \item Global minimum at $\vec{x^*} =$ (1, 1, 1, 1), $\vec{y^*} =$ (1,1,1,1,1,3,3,3), $f(\vec{x^*}) = -194$.
\end{itemize}

% \subsection{2.10-TP9 Problem}
\subsection{}
\begin{itemize}
    \item Dimension: $10$;
    \item Domain and search space: $0 \leq x_{i}$;
    \item Function: 
    \begin{equation}
    % f(\vec{x}) = -\left(\sum_{i=1}^{n-1}x_{i}x_{i+1} + \sum_{i=1}^{n-2}x_{i}x_{i+2} + x_{1}x_{9} + x_{1}x_{10} + x_{2}x_{10} + x_{1}x_{5} + x_{4}x_{7}\right);
    f(\vec{x}) = -\sum_{i=1}^{n-1}x_{i}x_{i+1} - \sum_{i=1}^{n-2}x_{i}x_{i+2} - x_{1}x_{9} - x_{1}x_{10} - x_{2}x_{10} - x_{1}x_{5} - x_{4}x_{7};
    \end{equation}
    \item Constraints: 
    \begin{equation}
    c_{1}(\vec{x}) = \sum_{i=1}^{n}x_{i} - 1 = 0
    \end{equation}
    \item Global minimum at $\vec{x^*} =$ (0, 0, 0, 0.25, 0.25, 0.25, 0.25, 0, 0, 0), $f(\vec{x^*}) = -0.375$.
\end{itemize}

% \subsection{3.4-TP3 Problem}
\subsection{}
\begin{itemize}
    \item Dimension: $6$;
    \item Domain and search space: $l_{i} \leq x_{i} \leq u_{i}$, with $l=(0,0,1,0,1,0)$, $u=(-,-,5,6,5,10)$;
    \item Function: 
    \begin{equation}
    f(\vec{x}) = -25(x_{1} - 2)^2 - (x_{2} - 2)^2 - (x_{3} - 1)^2 -(x_{4} - 4)^2 - (x_{5} - 1)^2 - (x_{6} - 4)^2;
    \end{equation}
    \item Constraints: 
    \begin{equation}\begin{cases}
    c_{1}(\vec{x}) = 4 - (x_{3} - 3)^2 - x_{4} \leq 0 \\
    c_{2}(\vec{x}) = 4 - (x_{5} - 3)^2 - x_{6} \leq 0 \\
    c_{3}(\vec{x}) = x_{1} - 3x_{2} -2 \leq 0 \\
    c_{4}(\vec{x}) = -x_{1} + x_{2} -2 \leq 0 \\
    c_{5}(\vec{x}) = x_{1} + x_{2} -6 \leq 0 \\
    c_{6}(\vec{x}) = 2 - x_{1} - x_{2} \leq 0 
    \end{cases}\end{equation}
    \item Global minimum at $\vec{x^*} =$ (5, 1, 5, 0, 5, 10), $f(\vec{x^*}) = -310$.
\end{itemize}

% \subsection{3.5-TP4 Problem}
\subsection{}
\begin{itemize}
    \item Dimension: $3$;
    \item Domain and search space: $l_{i} \leq x_{i} \leq u_{i}$, with $l=(0,0,0)$, $u=(2,-,3)$;
    \item Function: 
    \begin{equation}
    f(\vec{x}) = -2x_{1} + x_{2} - x_{3};
    \end{equation}
    \item Constraints: 
    \begin{equation}\begin{cases}
    c_{1}(\vec{x}) = -\vec{x}^T\vec{A}^T\vec{A}\vec{x} +2\vec{y}^T\vec{A}\vec{y} - \lVert \vec{y} \rVert^2 +0.25 \lVert \vec{b}-\vec{z} \rVert^2 \leq 0 \\
    c_{2}(\vec{x}) = x_{1} + x_{2} + x_{3} -4 \leq 0 \\
    c_{3}(\vec{x}) = 3x_{2} + x_{3} - 6 \leq 0
    \end{cases}\end{equation}
    \begin{equation*}
    \text{with} \quad \vec{A} = \begin{pmatrix} 0 & 0 & 1 \\ 0 & -1 & 0 \\ -2 & 1 & -1 \end{pmatrix}, \quad \vec{b}=(3,0,-4)^T
    \end{equation*}
    \begin{equation*}
    \vec{y} = (1.5,-0.5,-5)^T, \quad \vec{z}=(0,-1,-6)^T
    \end{equation*}
    \item Global minimum at $\vec{x^*} =$ (0.5, 0, 3), $f(\vec{x^*}) = -4$.
\end{itemize}

% \subsection{4.10-TP8 Problem}
\subsection{}
\begin{itemize}
    \item Dimension: $2$;
    \item Domain and search space: $l_{i} \leq x_{i} \leq u_{i}$, with $l=(0,0)$, $u=(2,3)$;
    \item Function: 
    \begin{equation}
    f(\vec{x}) = -12x_{1} - 7x_{2} + x_{2}^2;
    \end{equation}
    \item Constraints: 
    \begin{equation}
    c_{1}(\vec{x}) = -2x_{1}^4 + 2 - x_{2} = 0 
    \end{equation}
    \item Global minimum at $\vec{x^*} =$ (0.7175, 1.47), $f(\vec{x^*}) = -16.73889$.
\end{itemize}

% \subsection{4.10-TP9 Problem}
\subsection{}
\begin{itemize}
    \item Dimension: $2$;
    \item Domain and search space: $l_{i} \leq x_{i} \leq u_{i}$, with $l=(0,0)$, $u=(3,4)$;
    \item Function: 
    \begin{equation}
    f(\vec{x}) = -x_{1} - x_{2};
    \end{equation}
    \item Constraints: 
    \begin{equation}\begin{cases}
    c_{1}(\vec{x}) = x_{2} - 2 -2x_{1}^4 + 8x_{1}^3 - 8x_{1}^2 \leq 0 \\
    c_{2}(\vec{x}) = x_{2} - 4x_{1}^4 +32x_{1}^3 -88x_{1}^2 +96x_{1} -36 \leq 0
    \end{cases}\end{equation}
    \item Global minimum at $\vec{x^*} =$ (2.3295, 3.17846), $f(\vec{x^*}) = -5.50796$.
\end{itemize}

% \subsection{7.2.5-TP5 Problem}
\subsection{}
\begin{itemize}
    \item Dimension: $5$;
    \item Domain and search space: $l_{i} \leq x_{i} \leq u_{i}$, with $l=(78,33,27,27,27)$, $u=(102,45,45,45,45)$;
    \item Function: 
    \begin{equation}
    f(\vec{x}) = 5.3578x_{3}^2 + 0.8357x_{1}x_{5} + 37.2392x_{1};
    \end{equation}
    \item Constraints: 
    \begin{equation}\begin{cases}
    c_{1}(\vec{x}) = 0.00002584x_{3}x_{5} - 0.00006663x_{2}x_{5} - 0.0000734x_{1}x_{4} -1 \leq 0 \\
    c_{2}(\vec{x}) = 0.000853007x_{2}x_{5} + 0.00009395x_{1}x_{4} - 0.00033085x_{3}x_{5} -1 \leq 0 \\
    c_{3}(\vec{x}) = 1330.3294x_{2}^{-1}x_{5}^{-1} - 0.42x_{1}x_{5}^{-1} - 0.30586x_{2}^{-1}x_{3}^2x_{5}^{-1} -1 \leq 0 \\
    c_{4}(\vec{x}) = 0.00024186x_{2}x_{5} + 0.00010159x_{1}x_{2} + 0.00007379x_{3}^2 -1 \leq 0 \\
    c_{5}(\vec{x}) = 2275.1327x_{3}^{-1}x_{5}^{-1} - 0.2668x_{1}x_{5}^{-1} - 0.40584x_{4}x_{5}^{-1} -1 \leq 0 \\
    c_{6}(\vec{x}) = 0.00029955x_{3}x_{5} + 0.00007992x_{1}x_{3} + 0.00012157x_{3}x_{4} -1 \leq 0
    \end{cases}\end{equation}
    \item Global minimum at $\vec{x^*} =$ (78, 33, 29.998, 45, 36.7673), $f(\vec{x^*}) = 10122.696$.
\end{itemize}

% \subsection{7.2.6-TP6 Problem}
\subsection{}
\begin{itemize}
    \item Dimension: $3$;
    \item Domain and search space: $1 \leq x_{i} \leq 100$;
    \item Function: 
    \begin{equation}
    f(\vec{x}) = 0.5x_{1}x_{2}^{-1} - x_{1} - 5x_{2}^{-1};
    \end{equation}
    \item Constraints: 
    \begin{equation}
    c_{1}(\vec{x}) = 0.01x_{2}x_{3}^{-1} + 0.01x_{1} +0.0005x_{1}x_{3} -1 \leq 0
    \end{equation}
    \item Global minimum at $\vec{x^*} =$ (88.2890, 7.7737, 1.3120), $f(\vec{x^*}) = -83.254$.
\end{itemize}

% \subsection{7.2.7-TP7 Problem}
\subsection{}
\begin{itemize}
    \item Dimension: $4$;
    \item Domain and search space: $0.1 \leq x_{i} \leq 10$;
    \item Function: 
    \begin{equation}
    f(\vec{x}) = -x_{1} + 0.4x_{1}^{0.67}x_{3}^{-0.67};
    \end{equation}
    \item Constraints: 
    \begin{equation}\begin{cases}
    c_{1}(\vec{x}) = 0.05882x_{3}x_{4} + 0.1x_{1} -1 \leq 0 \\
    c_{2}(\vec{x}) = 4x_{2}x_{4}^{-1} + 2x_{2}^{-0.71}x_{4}^{-1} + 0.05882x_{2}^{-1.3}x_{3} -1 \leq 0
    \end{cases}\end{equation}
    \item Global minimum at $\vec{x^*} =$ (8.1267, 0.6154, 0.5650, 5.6368), $f(\vec{x^*}) = -5.7398$.
\end{itemize}

% \subsection{7.2.8-TP8 Problem}
\subsection{}
\begin{itemize}
    \item Dimension: $8$;
    \item Domain and search space: $0.01 \leq x_{i} \leq 10$;
    \item Function: 
    \begin{equation}
    f(\vec{x}) = -x_{1} - x_{5} + 0.4x_{1}^{0.67}x_{3}^{-0.67} + 0.4x_{5}^{0.67}x_{7}^{-0.67};
    \end{equation}
    \item Constraints: 
    \begin{equation}\begin{cases}
    c_{1}(\vec{x}) = 0.05882x_{3}x_{4} + 0.1x_{1} -1 \leq 0 \\
    c_{2}(\vec{x}) = 0.05882x_{7}x_{8} + 0.1x_{1} + 0.1x_{5} -1 \leq 0 \\
    c_{3}(\vec{x}) = 4x_{2}x_{4}^{-1} + 2x_{2}^{-0.71}x_{4}^{-1} + 0.05882x_{2}^{-1.3}x_{3} -1 \leq 0 \\
    c_{4}(\vec{x}) = 4x_{6}x_{8}^{-1} + 2x_{6}^{-0.71}x_{8}^{-1} + 0.05882x_{6}^{-1.3}x_{7} -1 \leq 0
    \end{cases}\end{equation}
    \item Global minimum at $\vec{x^*} =$ (6.4225, 0.6686, 1.0239, 5.9399, 2.2673, 0.5960, 0.4029, 5.5288), $f(\vec{x^*}) =$ -6.0482.
\end{itemize}

% \subsection{8.2.7-TP7 Problem}
\subsection{}
\begin{itemize}
    \item Dimension: $5$;
    \item Domain and search space: $-5 \leq x_{i} \leq 5$;
    \item Function: 
    \begin{equation}
    f(\vec{x}) = (x_{1}-1)^2 + (x_{1}-x_{2})^2 + (x_{2}-x_{3})^3 + (x_{3}-x_{4})^4 + (x_{4}-x_{5})^4;
    \end{equation}
    \item Constraints: 
    \begin{equation}\begin{cases}
    c_{1}(\vec{x}) = x_{1} + x_{2}^2 + x_{3}^3 -3\sqrt{2} -2 = 0 \\
    c_{2}(\vec{x}) = x_{2} - x_{3}^2 + x_{4} -2\sqrt{2} +2 = 0 \\
    c_{3}(\vec{x}) = x_{1}x_{5} - 2 = 0
    \end{cases}\end{equation}
    \item Global minimum at $\vec{x^*} =$ (1.1166, 1.2204, 1.5378, 1.9728, 1.7911), $f(\vec{x^*}) =$ 0.0293.
\end{itemize}

\end{document}